\documentclass[12pt]{amsart}
\usepackage{amsthm, amssymb, dsfont}
\usepackage{hyperref}
\usepackage[all]{xy}
\usepackage[active]{srcltx}
\usepackage[short,nodayofweek]{datetime}
\usepackage{xcolor, enumitem}

\usepackage{tikz}
\usetikzlibrary{arrows}

\definecolor{my-blue}{cmyk}{1,0.6,0,0}
\definecolor{my-green}{cmyk}{0.8,0,1,0.5}

\hypersetup{
colorlinks=true, 
linkcolor=my-blue, 
citecolor=my-green,
urlcolor=black
}

\makeatletter
\@namedef{subjclassname@2020}{%
  \textup{2020} Mathematics Subject Classification}
\makeatother

\theoremstyle{plain}
\newtheorem{thm}{Theorem}[section]
\newtheorem{cor}[thm]{Corollary}
\newtheorem{lem}[thm]{Lemma}
\newtheorem{prop}[thm]{Proposition}

\theoremstyle{definition}
\newtheorem{defn}[thm]{Definition}
\newtheorem{exmp}[thm]{Example}
\newtheorem{rem}[thm]{Remark}

\newtheorem{algo}[thm]{Algorithm}
\newtheorem*{questions}{Questions}

\numberwithin{equation}{section}

\newcommand{\Ga}{{\mathbb G}_a}
\newcommand\NN{{\mathbb N}}
\newcommand\ZZ{{\mathbb Z}}
\newcommand{\Fq}{{{\mathbb F}_{\!q}}}
\newcommand{\cD}{\mathcal{D}}
\newcommand{\cP}{\mathcal{P}}
\newcommand{\cF}{\mathcal{F}}

\newcommand{\mot}{\mathsf{M}}  
\newcommand{\dumot}{\mathfrak{M}}  
\newcommand{\grp}{\mathrm{grp}}

\newcommand{\fop}{\pi}  
\newcommand{\sop}{\rho}  
\newcommand{\Kpolys}{K\{\fop,\sop\}}
\newcommand{\lm}{\mathrm{lm}}  
\newcommand{\lc}{\mathrm{lc}}  
\newcommand{\tor}[1]{#1\text{-tor}}
\newcommand{\Jtop}{J_{\rm{top}}}
\newcommand{\Jlow}{J_{\rm{low}}}
\newcommand{\Btop}{B_{\rm{top}}}
\newcommand{\Blow}{B_{\rm{low}}}

\newcommand{\one}{\mathds{1}}
\newcommand{\bigmid}{\, \big| \,}
\newcommand{\dk}{\check{\kappa}} 
\newcommand{\bkappa}{\bar{\kappa}} 
\newcommand{\bdk}{\bar{\dk}} 
\newcommand{\be}{\bar{e}} 
\newcommand{\ce}{\check{e}} 

\newcommand{\sfn}{\mathsf{n}}
\newcommand{\sfm}{\mathsf{m}}

\newcommand{\Wgen}{W_{\rm{gen}}}
\newcommand{\bWgen}{\bar{W}_{\rm{gen}}}

\newcommand{\Wind}{W_{\rm{ind}}}
\newcommand{\bWind}{\bar{W}_{\rm{ind}}}

\newcommand{\rN}{\mathrm{N}}

\newcommand{\pr}{\mathrm{pr}}
\newcommand{\inj}{\mathrm{in}}
\newcommand{\dphi}{d\phi}
\newcommand{\id}{\mathrm{id}}
\newcommand{\rk}{\mathrm{rk}}
\newcommand{\diag}{\mathrm{diag}}

\renewcommand{\sp}[1]{\{#1\}} 
\newcommand{\sls}[1]{(\!\{#1\}\!)}

\newcommand{\partdef}[1]{ \left\{ \begin{array}{ll} #1 \end{array} \right. }

\newcommand{\gen}[1]{\langle #1 \rangle}
\newcommand{\Dgen}[1]{{}_\cD\gen{#1}}

\newcommand{\Ktgen}[1]{{}_{K[t]}\gen{#1}}
\newcommand{\Ksopgen}[1]{{}_{K\sp{\sop}}\gen{#1}}
\newcommand{\svect}[2]{\left( \begin{matrix} {#1}_{1}\\ \vdots \\ {#1}_{#2}\end{matrix}\right)}
\newcommand{\zvect}[2]{\left( \begin{matrix} {#1}_{1}& {#1}_{2} & \cdots & {#1}_{#2}\end{matrix}\right)}

\DeclareMathOperator{\Hom}{Hom}
\DeclareMathOperator{\End}{End}
\DeclareMathOperator{\Mat}{Mat}
\DeclareMathOperator{\GL}{GL}
\DeclareMathOperator{\Lie}{Lie}
\DeclareMathOperator{\Mon}{Mon}
\DeclareMathOperator{\NF}{NF}

\newcommand{\grid}[3]{
    \draw[axis] (-1,0)  -- (4,0) node [below right] {$#1$};
    \draw[axis] (0,-1) -- (0,4) node [above left] {$#2$};
    \draw[semithick] (1,-0.1) node [below] {$1$} -- (1,0.1); 
    \draw[semithick] (2,-0.1) node [below] {$2$} -- (2,0.1); 
    \draw[semithick] (3,-0.1) node [below] {$3$} -- (3,0.1); 

    \draw[semithick] (-0.1,1) node [left] {$1$} -- (0.1,1); 
    \draw[semithick] (-0.1,2) node [left] {$2$} -- (0.1,2); 
    \draw[semithick] (-0.1,3) node [left] {$3$}-- (0.1,3); 
    \draw (2.5, 4.5) node {$\kappa_{#3}$};
}
\newcommand{\linetau}[2]{
    \draw [color=red, thick] (#1,#2) -- (3.8,#2);
}
\newcommand{\linet}[2]{
    \draw [color=red, thick] (#1,#2) -- (#1,3.8);
}

\newcommand{\cone}[2]{
	\linet{#1}{#2}
	\linetau{#1}{#2}
    \draw [fill=red, opacity=0.15] (#1,#2) rectangle (3.8,3.8);
}
\newcommand{\point}[4]{
    \draw [fill, color=red] (#3,#4) circle [radius=.08] node [#2] {$\textcolor{red}{#1}$};
}
\newcommand{\myscale}{0.7}

\newcommand{\twosheets}[2]{
\hspace*{2cm}
 \begin{tikzpicture}[
    scale=\myscale,
    axis/.style={thin, ->, >=stealth'},
    every node/.style={color=black},
    ]
    #1
\end{tikzpicture} \hfill
\begin{tikzpicture}[
    scale=\myscale,
    axis/.style={thin, ->, >=stealth'},
    every node/.style={color=black},
    ]
	#2
\end{tikzpicture} 
\hspace*{2cm}}

\begin{document}
\baselineskip=17pt

\title[Determining $t$-motives]{Determining $t$-motives and dual $t$-motives in Anderson's theory}
\author{Andreas Maurischat}
\address{Faculty of Computer Science\\ RWTH Aachen University\\ 52072 Aachen, Germany}
\email{maurischat@combi.rwth-aachen.de}

\newdate{date}{22}{01}{2026}
\date{\displaydate{date}}

%


\begin{abstract}
Anderson t-modules are analogs of abelian varieties in positive characteristic. Associated to such a t-module, there are its t-motive and its dual t-motive. 
When dealing with these objects, several questions occur which one would like to solve algorithmically. For example, for a given t-module one would like to decide whether its t-motive is indeed finitely generated free, and determine a basis. Reversely, for a given object in the category of t-motives one would like to decide whether it is the t-motive associated to a t-module, and determine that t-module.

In this article, we positively answer such questions by providing the corresponding algorithms.

As it turned out, the main part of all these algorithms stems from a single algorithm in non-commutative algebra, and hence the first part of this article doesn't deal with Anderson's objects at all, but are results on finitely generated modules over skew polynomial rings.

\end{abstract}

\maketitle

\setcounter{tocdepth}{1}
\tableofcontents

\section*{Introduction}

In function field arithmetic, key objects to study are Anderson $t$-modules, $t$-motives, and dual $t$-motives -- the function field analogs in positive characteristic of abelian varieties and motives.
When dealing with these objects, several questions occur which one would like to solve algorithmically.
For explaining these questions, we have to introduce some terminology.

Let $K$ be a field containing the finite field $\Fq$ of $q$ elements. Roughly speaking, an \emph{Anderson $t$-module} over $K$ of dimension $d$ is a certain pair $(E,\phi)$ where $E$ is an algebraic group over $K$ which is isomorphic to $\Ga^d$ -- the $d$-fold product of the additive group --, equipped with a compatible $\Fq$-action, and 
\[\phi:\Fq[t]\to \End_{\grp,\Fq}(E), a\mapsto \phi_a\]
is a certain $\Fq$-algebra homomorphism into $\Fq$-linear endomorphisms of $E$.\footnote{See Section \ref{sec:t-modules} for the precise definition.}

There are two kinds of motives attached to Anderson $t$-modules.
The first one -- called \emph{$t$-motive} -- was defined by Anderson in his seminal paper \cite{ga:tm}. Among others, it was used ibid.~to provide a new criterion for uniformizability of $t$-modules.
The other kind -- usually called \emph{dual $t$-motive} -- was defined later in \cite{ga-wb-mp:darasgvpc}, and is the base for the ABP-criterion \cite[Thm.~3.1.1]{ga-wb-mp:darasgvpc} and Papanikolas' Theorem \cite[Thm.~5.2.2]{mp:tdadmaicl} which opened a new way to study transcendence questions in positive characteristic.

As the notion ``dual $t$-motive'' might cause confusion with ``dual of a $t$-motive'', we will rather use the terminology \emph{$t$-comotive} as in \cite{qg-am:pamfrfe}.

Both $t$-motives and $t$-comotives are $K[t]$-modules, and the $t$-module $E$ is said to be \emph{abelian} if its $t$-motive $\mot(E)$ is a finitely generated $K[t]$-module, and the $t$-module $E$ is said to be \emph{$t$-finite} (or for consistency \emph{coabelian}) if its $t$-comotive $\dumot(E)$ is a finitely generated $K[t]$-module.

In \cite[Theorem A]{am:aefam}, we showed that these two notion are indeed equivalent.

Anderson already showed that for an abelian (resp.~coabelian) $t$-module, its $t$-motive (resp.~its $t$-comotive) is even a free $K[t]$-module. So for computational purposes one is also interested in obtaining a $K[t]$-basis of the $t$-motive (resp.~the $t$-comotive), as well as a description of the Frobenius twist action induced by the Frobenius action on the additive group $\Ga$.

Although \cite[Theorem A]{am:aefam} provides a computable criterion to check whether a $t$-module is abelian/coabelian, it does not solve these two tasks.

Coming from the different direction, one can start with an object from the category of $t$-motives which is defined to consist of finitely generated free $K[t]$-modules with a semilinear action of the Frobenius twist operator 
\[ \tau:K[t]\to K[t],\sum x_it^i\mapsto \sum x_i^qt^i,\] and one is faced with the question whether this $t$-motive is the $t$-motive of some Anderson $t$-module. And if it is, how can the $t$-action of the Anderson $t$-module be described.

The same questions occur for the $t$-comotive. Here, one has a semilinear action of the inverse Frobenius twist operator $\sigma=\tau^{-1}$ and has the extra assumption that $K$ is perfect for guaranteeing the well-definedness of $\sigma$.

When considering tensor products of Anderson $t$-modules, one is even forced to walk both directions, since the tensor product of Anderson $t$-modules is defined via the $K[t]$-tensor product of the associated $t$-motives or of the associated $t$-comotives.\footnote{We would like to remind the reader that for more general Anderson $A$-modules, one would have to take care which of these two constructions of the tensor product is used, since the results are only isogeneous, and not neccessarily isomorphic (see e.g.~\cite[Cor.~6.4]{qg-am:pamfrfe}). For Anderson $t$-modules, however, these two constructions give indeed isomorphic $t$-modules.
}

We summarize the mentioned tasks as the following

\begin{questions} \ 
\begin{enumerate}[label=\arabic*., ref=\arabic*]
\item \label{item:basis-of-t-motive} Given explicitly\footnote{What we mean by ``Given explicitly'' in these questions will become clear in the later sections.} an Anderson $t$-module $(E,\phi)$, how can we
\begin{enumerate}
\item check whether $E$ is abelian,
\item compute a $K[t]$-basis of $\mot(E)$ (in case that $E$ is abelian), and
\item describe the Frobenius twist action with respect to this $K[t]$-basis?
\end{enumerate}
\item \label{item:t-module-to-t-motive} Given explicitly an effective Anderson $t$-motive $\mot$, how can we
\begin{enumerate}
\item decide whether $\mot$ is (isomorphic to) the $t$-motive associated to some Anderson $t$-module $(E,\phi)$,
\item compute an isomorphism $E\to \Ga^d$ (if such an $E$ exists), and
\item describe $\phi$ with respect to this isomorphism?
\end{enumerate}
\item \label{item:basis-of-dual-t-motive} Given explicitly an Anderson $t$-module $(E,\phi)$ over a perfect field $K$, how can we
\begin{enumerate}
\item check whether $E$ is coabelian,
\item compute a $K[t]$-basis of $\dumot(E)$ (in case that $E$ is coabelian), and
\item describe the inverse Frobenius twist action with respect to this $K[t]$-basis?
\end{enumerate}
\item \label{item:t-module-to-dual-t-motive} Given explicitly an effective Anderson $t$-comotive $\dumot$, how can we
\begin{enumerate}
\item decide whether $\dumot$ is (isomorphic to) the $t$-comotive associated to some Anderson $t$-module $(E,\phi)$,
\item compute an isomorphism $E\to \Ga^d$ (if such an $E$ exists), and
\item describe $\phi$ with respect to this isomorphism?
\end{enumerate}
\end{enumerate}
\end{questions}

In this article, we positively answer all these questions (see Theorems \ref{thm:criterion-for-abelianess}, \ref{thm:criterion-for-t-module-to-motive}, \ref{thm:criterion-for-t-finiteness} and \ref{thm:criterion-for-t-module-to-dual-motive}).

We roughly sketch the main ideas to the answers.

The $t$-motive $\mot$ can be considered as a module over the skew polynomial ring in two variables $K\sp{\tau,t}$. By hypothesis it is free of finite rank over the skew polynomial ring $K\sp{\tau}$ in Question \ref{item:basis-of-t-motive}, and free of finite rank over the polynomial ring $K[t]$ in Question \ref{item:t-module-to-t-motive}, in both cases with given basis and given action on the basis of the other variable.

The answers to the questions are then closely related to finding out whether the $t$-motive $\mot$ is free and finitely generated over the other one-variable polynomial ring ($K[t]$ in the first case, and $K\sp{\tau}$ in the second), and if it is, to computing a basis over that polynomial ring, and giving the $\tau$-action (resp.~the $t$-action) with respect to that basis.

The same holds for the questions on the $t$-comotive with $\tau$ replaced by $\sigma$.

These considerations motivated the setting in Section \ref{sec:notation}, where we consider finitely generated (left) modules $M$ over a 
skew polynomial ring $\cD=\Kpolys$ in two variables $\fop$ and~$\sop$. As $\cD$ is left-Noetherian such an $M$ can be given by finitely many generators and finitely many relations, i.e., we have an isomorphism 
\[ M\cong \cF/\Dgen{p_1,\ldots, p_d} \]
giving $M$ as the quotient of a free $\cD$-module $\cF$ of finite rank by the submodule generated by the elements $p_1,\ldots, p_d$.

The key to the answers is to use an adaption of Janet's algorithm to linear difference equations as given in \cite{dr:faepde} which we recall in Section \ref{sec:janet-basis-and-algorithm}.
Janet’s algorithm turns the set $\{p_1,\ldots,p_d\}$ of generators for the relations into a new set
of generators $\{b_1,\ldots, b_r \}$ a so called \emph{Janet basis} such that unique normal forms for the representatives of each residue class in $M$ can be computed.

The output of Janet's algorithm and the normal forms depend on a chosen monomial order on $\Kpolys$ and on $\cF$. 
For our application, the lexicographical order on $\Kpolys$ with $\sop \prec \fop$, and its extension to $\cF$ via a position-over-term order is exactly what we need.
Namely with this chosen monomial order, on can readily read off from the output of the algorithm (see Theorem \ref{thm:criterion-for-finiteness} and Theorem \ref{thm:relations-and-fop-action})
\begin{itemize}
\item whether $M$ is finitely generated free as $K\sp{\sop}$-module, and what its rank is, and in that case
\item a finite $K\sp{\sop}$-generating set, and the $\fop$-action on this set, as well as
\item the $K\sp{\sop}$-relations among this generating set.
\end{itemize}

The final step to obtain a $K\sp{\sop}$-basis and the $\fop$-action on that basis, is obtained by applying the (non-commutative version of) the elementary divisor algorithm on that data (see Remark \ref{rem:finding-a-basis}).

For answering the questions stated above, we apply this procedure to the various settings, and only have to do small extra computations in order to rewrite data for the Anderson $t$-module to corresponding data for the $t$-motive or the $t$-comotive, and vice versa.

The manuscript is structured as follows.
In Section \ref{sec:notation}, we introduce the basic algebraic notation of our setting. The definition of Janet bases and the description of the Janet algorithm simplified to our situation are given in Section \ref{sec:janet-basis-and-algorithm}.
Section \ref{sec:finite-generation-and-relations} is dedicated to the main theorems for extracting the relevant data from a Janet basis.
These first three sections provide results in non-commutative algebra, and don't deal with Anderson's objects at all.

In Section \ref{sec:t-modules}, we recall the objects of Anderson's theory and their interrelation. Then we answer Questions \ref{item:basis-of-t-motive} and \ref{item:t-module-to-t-motive} in Sections \ref{sec:basis} and \ref{sec:reverse-direction}, respectively. Finally, Questions \ref{item:basis-of-dual-t-motive} and \ref{item:t-module-to-dual-t-motive} on the $t$-comotive are answered in Section \ref{sec:dual-t-motives}.

The answer to the first question is illustrated by several examples in Section \ref{sec:examples}.

\section{Notation and basics}\label{sec:notation}

Throughout this article, $K$ is a field, and $\Kpolys$ is the skew polynomial ring in two variables $\fop$ and $\sop$, such that for all $x\in K$,
\[  \fop\cdot \sop = \sop\cdot \fop,\quad \fop\cdot x = \gamma_\fop(x)\cdot \fop,\quad \sop\cdot x = \gamma_\sop(x)\cdot \sop, \]
where $\gamma_\fop,\gamma_\sop:K\to K$ are two commuting endomorphisms of fields. 

Every element in $\Kpolys$ can therefore uniquely be written as a sum of monomial terms $x\fop^k\sop^j$ where $x\in K$, $k,j\geq 0$.

At certain places, we will additionally assume that $\gamma_\sop$ is an automorphism, but $\gamma_\fop$ is always allowed to not be surjective.

The skew polynomial ring $\Kpolys$ is a left-Noetherian ring, and so are its subrings $K\sp{\fop}$ and $K\sp{\sop}$. These two subrings are even left principal ideal domains, i.e.,~every left ideal is generated by one element. This is due to the fact that one can always perform a right-division with remainder in these two subrings, and a right Euclidean algorithm to obtain a greatest common right-hand divisor of two elements (see \cite[Section 2]{oo:tncp}).

Throughout this article, we will consider left modules, but will usually omit the word~``left''.

We start with an observation:
\begin{prop}\label{prop:observations-on-M}
Let $M$ be a $\Kpolys$-module, and denote by $M_{\tor{\fop}}$ and $M_{\tor{\sop}}$ the $K\sp{\fop}$-torsion subspace, and the $K\sp{\sop}$-torsion subspace respectively, i.e., the $K$-vector spaces
\begin{align*}
M_{\tor{\fop}} &:= \left\{ f\in M \mid \exists\, 0\ne p\in K\sp{\fop}: p\cdot f=0\right\},\\ 
M_{\tor{\sop}} &:= \left\{ f\in M \mid \exists\, 0\ne p\in K\sp{\sop}: p\cdot f=0\right\}. 
\end{align*}
Then
\begin{enumerate}
\item \label{item:torsion-is-submodule} $M_{\tor{\fop}}$ and $M_{\tor{\sop}}$ are $\Kpolys$-submodules.
\item \label{item:inclusion-of-torsion} If $M$ is finitely generated as a $K\sp{\sop}$-module, then $M_{\tor{\sop}}\subseteq M_{\tor{\fop}}$.
\item \label{item:freeness} If $M$ is finitely generated as a $K\sp{\sop}$-module, torsionfree as $K\sp{\fop}$-module, and $\gamma_\sop$ is an automorphism, then $M$ is a free $K\sp{\sop}$-module.
\end{enumerate}
\end{prop}

\begin{rem}
Part \eqref{item:freeness} of the proposition is what we will need in the later parts. Special cases of that statement are already given in \cite[Lemma 1.4.5]{ga:tm} and \cite[Proposition 4.3.2]{ga-wb-mp:darasgvpc}. Actually, for the applications to $t$-motives and $t$-comotives these special cases are sufficient. 
\end{rem}

\begin{proof}[Proof of Proposition \ref{prop:observations-on-M}] \
\begin{enumerate}
\item For $f\in M_{\tor{\fop}}$, let $0\ne p\in K\sp{\fop}$ be such that $p\cdot f=0$. Then
\begin{align*}
0 &= \fop\cdot p\cdot f =\gamma_\fop(p)\cdot ( \fop \cdot f),\quad \text{and}\\
0 &= \sop\cdot p\cdot f =\gamma_\sop(p)\cdot ( \sop \cdot f),
\end{align*}
where by abuse of notation, $\gamma_\fop(p)$ and $\gamma_\sop(p)$ denote the polynomials obtained by applying $\gamma_\fop$ and $\gamma_\sop$, respectively, to the coefficients of $p$. Since, $0\ne\gamma_\fop(p)\in K\sp{\fop}$ and $0\ne\gamma_\sop(p)\in K\sp{\fop}$, this shows that $\fop\cdot f,\sop\cdot f\in M_{\tor{\fop}}$. Hence $M_{\tor{\fop}}$ is a $\Kpolys$-module. The proof for $M_{\tor{\sop}}$ is the same.
\item Assume that $M$ is finitely generated as $K\sp{\sop}$-module. Since $K\sp{\sop}$ is left-Noetherian, all $K\sp{\sop}$-submodules of $M$ are also finitely generated, and in particular $M_{\tor{\sop}}$.

As for $f\in M_{\tor{\sop}}$, one has $\dim_K( K\sp{\sop}f)<\infty$, finite generation of $M_{\tor{\sop}}$ implies that $\dim_K( M_{\tor{\sop}} )<\infty$.

For any $0\ne f\in M_{\tor{\sop}} $, we have $\Kpolys f \subseteq M_{\tor{\sop}}$ by Part \eqref{item:torsion-is-submodule}, and in particular, $K\sp{\fop}f\subseteq M_{\tor{\sop}}$. Hence $\dim_K(K\sp{\fop}f)\leq \dim_K(M_{\tor{\sop}})<\infty$, which implies that $f\in M_{\tor{\fop}}$.
\item By Part \eqref{item:inclusion-of-torsion}, we readily see that $M_{\tor{\sop}}=0$, i.e., that $M$ is torsionfree as $K\sp{\sop}$-module.
When $\gamma_\sop$ is an automorphism, the ring $K\sp{\sop}$ is also a right principal ideal domain, and the structure theorem on finitely generated modules over (non-commutative) principal ideal domains (see \cite[Theorem 18]{nj:tr}) ensures that $M$ is free.
\end{enumerate}
\end{proof}


\section{Janet basis and algorithm}\label{sec:janet-basis-and-algorithm}

This section is mainly an excerpt of \cite[Section 2.1]{dr:faepde}, adapted to our situation.
We denote by $\cD:=\Kpolys$ the skew polynomial ring defined in Section \ref{sec:notation}.

We first sketch the rough idea of Janet bases and of Janet's algorithm in algebraic terms.

The starting point is a finitely generated (left) module $M$ over the ring $\cD$.
As $\cD$ is left-Noetherian, $M$ is given by a finite set of generators and a finite set of relations.
Hence,  we have an isomorphism
\begin{equation}\label{eq:presentation-of-M}
\cF/\Dgen{g_1,\ldots, g_s} \xrightarrow{\cong} M
\end{equation} 
of $\cD$-modules, where $\cF=\bigoplus_{i=1}^d \cD\kappa_i$ is a free $\cD$-module with some basis $\{\kappa_1,\ldots, \kappa_d\}$, and
$\Dgen{g_1,\ldots, g_s}$ denotes the submodule generated by certain elements $g_1,\ldots, g_s\in \cF$.

We usually write $\bkappa_i$ for the image of $\kappa_i$ in $M$.

Janet's algorithm turns this set $\{g_1,\ldots, g_s\}$ of generators for the relations into a new set of generators $\{b_1,\ldots,b_r\}$ a so called \emph{Janet basis} such that
\begin{itemize}
\item a (unique) normal form for the representatives of each residue class in $M$ is defined and 
\item these normal forms can be computed effectively.
\end{itemize}

Originally, Janet formulated this algorithm in the case that $\cD$ is a ring of differential operators (see \cite{mj:lsedp}), but it was generalized to rings of difference operators by Robertz (see \cite[Section 2.1.3]{dr:faepde}). For simplicity, in our presentation here, we restrict to the special case of $\cD=\Kpolys$, and we fix a monomial order on $\cD$ as well as a specific order on monomials in $\cF$.\footnote{Robertz uses automorphisms, but everything works for any endomorphisms, since only right division with remainders are used.} 

\subsection{Monomial orders and normal forms}
We use the following notation:
\begin{itemize}
\item $\Mon(\{\fop,\sop\})=\{ \fop^k\sop^j \mid k,j\geq 0\}$ is the set of \emph{monomials} in $\cD=\Kpolys$, i.e.,~monomials in the variables $\fop$ and $\sop$, and accordingly,
\begin{align*}
\Mon(\{\fop\}) &= \{ \fop^k \mid k\geq 0\},\\
\Mon(\{\sop\}) &= \{ \sop^j \mid j\geq 0\},\\
\Mon(\emptyset) &= \{ 1 \}.
\end{align*}
\item In $\Kpolys$, we use the lexicographical order on monomials $\fop^k\sop^j$ where $\sop\prec \fop$, i.e., 
\[ \fop^{k_1}\sop^{j_1}\prec \fop^{k_2}\sop^{j_2} \Longleftrightarrow k_1<k_2 \text{ or } (k_1=k_2 \text{ and } j_1<j_2),\]
e.g.~$\sop\prec \fop \prec \fop^2 \prec \fop^2 \sop\prec \fop^3 \prec \fop^3 \sop$.
\item $\cF$ is a free left $\cD$-module of dimension $d$ with a fixed basis $\{\kappa_1,\ldots, \kappa_d\}$.
\item Monomials in $\cF$ are $m\cdot \kappa_i$ with $m$ a monomial in $\Kpolys$ and $\kappa_i$ being one of the basis elements.
\item On monomials in $\cF$, we use the \emph{position over term order} defined by\[
m_1\kappa_{i_1}\prec m_2\kappa_{i_2} \Longleftrightarrow i_1>i_2 \text{ or } (i_1=i_2 \text{ and } m_1\prec m_2 ) \] 
\item For $f\in \cF\setminus \{0\}$ its \emph{leading monomial}, denoted by $\lm(f)$, is the greatest monomial that occurs in the standard representation $f=\sum_{k,j,i} c_{k,j,i} \fop^k \sop^j\kappa_i$ of $f$ with non-zero coefficient. The corresponding coefficient $c_{k,j,i}$ will be called the \emph{leading coefficient} of $f$, denoted by $\lc(f)$.
\item For a monomial $m\kappa_i$ in $\cF$ and a subset $\mu\subset \{\fop,\sop\}$, we call the set of monomials
$\Mon(\mu)m\kappa_i = \{ \tilde{m}m\kappa_i \mid \tilde{m}\in \Mon(\mu) \}$, the \emph{$\mu$-cone} of $m\kappa_i$.
\item For $\{g_1,\ldots, g_s\}\subseteq \cF$, we denote by $\Dgen{g_1,\ldots,g_s}$ the $\cD$-submodule of $\cF$ generated by $g_1,\ldots, g_s$.
\end{itemize}

For defining a Janet basis for the submodule $\Dgen{g_1,\ldots,g_s}$, we have to consider finite sets $T=\{ (b_1,\mu_1),\ldots,(b_r,\mu_r)\}\subseteq \cF\times \cP(\{\fop,\sop\})$ where $\cP(\{\fop,\sop\})$ denotes the power set of $\{\fop,\sop\}$, i.e.,~the set of subsets of $\{\fop,\sop\}$.

\begin{defn}
For $T=\{ (b_1,\mu_1),\ldots,(b_r,\mu_r)\}\subseteq \cF\times \cP(\{\fop,\sop\})$, and $g\in \cF$, we define the \emph{reduction of $g$ modulo $T$} (or the \emph{normal form} of $g$ with respect to $T$), denoted $\NF(g,T)$, as the result of the following process:\\
Start by setting $h=g$.
If there is a monomial $m\kappa_i$ in $h$ and $(b,\mu)\in T$ such that
$m\kappa_i\in \Mon(\mu)\lm(b)$, find the largest such monomial $m\kappa_i$, and a corresponding $(b,\mu)\in T$. Let $\tilde{m}\in \Mon(\mu)$ be such that $m\kappa_i=\tilde{m}\cdot \lm(b)$, and $c\in K$ such that $c\cdot \lc(\tilde{m} \cdot b)$ equals the coefficient in $h$ of $m\kappa_i$. Then replace $h$ by
$h - c\cdot \tilde{m} \cdot b$.

Repeat this process until no such monomial $m$ in $h$ exists anymore.

The resulting $h$ is the \emph{normal form} of $g$ with respect to $T$: $h=\NF(g,T)$.
\end{defn}

\begin{rem}
Be aware, that in general, for a given monomial $m\kappa_i$, there might be several corresponding pairs $(b,\mu)\in T$, and the result of the reduction process might depend on the choice of the pair. However, this ambiguity will not cause any problems, as it doesn't occur for Janet bases.
\end{rem}

\begin{rem}
We will also use the normal form $\NF(g,H)$ for some set $H=\{b_1,\ldots, b_r\}\subseteq \cF$. In that, we mean the reduction with respect to $\{ (b_1,\mu_1),\ldots,(b_r,\mu_r)\}$ where all $\mu_i=\{\fop,\sop\}$.
\end{rem}

\begin{rem}
For the reduction process, the important information on such a set $T$ are the $\mu$-cones of $\lm(b)$ for the elements $(b,\mu)\in T$.
Namely, an element $h\in \cF$ can be reduced, if it has a monomial that lies in the $\mu$-cone of $\lm(b)$ for some $(b,\mu)\in T$.
We explain now how we depict this information (see Figure \ref{fig:graphic-explanation}):

\begin{figure}[b]
\twosheets{
	\grid{\fop}{\sop}{1}
	\point{p_2}{above right}{2}{0}
	\point{p_1}{above right}{1}{1}
	\point{p_3}{above right}{0}{2}
    \cone{2}{0}
    \linet{1}{1} 
    \linet{0}{2} 
}{
	\grid{\fop}{\sop}{2}
	\point{p_4}{above left}{3}{0}
    \cone{3}{0}
}
\caption{Graphical visualization for the set
\[  T= \{ (p_1,\{\sop\}), (p_2,\{\fop,\sop\}), (p_3,\{\sop\}), (p_4,\{\fop,\sop\})\}, \text{ where} \]
\begin{align*}
	p_1 &= -\underline{\fop \sop\kappa_1}+\fop\kappa_1 + \fop \kappa_2+\sop\kappa_2,\quad &
	p_3 &= -(\underline{\sop^2}-2\sop+1)\underline{\kappa_1}+ (\fop^2-\sop)\kappa_2,\quad \\
	p_2 &= \underline{\fop^2\kappa_1}+\sop\kappa_1-\fop\kappa_2-\kappa_2,\quad &
	p_4 &= \underline{\fop^3\kappa_2}+\fop^2\sop\kappa_2 -\sop^2\kappa_2.
\end{align*} 
(The underlined terms are the leading monomials.)\label{fig:graphic-explanation}
}
\end{figure}

For each basis vector $\kappa_i$, we draw a two dimensional coordinate system with axes labelled by $\fop$ and $\sop$, which we call \emph{the $\kappa_i$-sheet} in the following. Each monomial $\fop^l \sop^j\kappa_i\in \Mon(\cF)$ is then identified with the point $(l,j)$ in the $\kappa_i$-sheet.

Now, for $(b,\mu)\in T$, the $\mu$-cone of $\lm(b)$ is depicted as follows: We color the point corresponding to $\lm(b)$, and label it by $b$. Further, if $\sop\in \mu$ (resp.~ $\fop\in \mu$), we add a colored half-line parallel to the $\sop$-axis (resp.~the $\fop$-axis) starting at this point. If $\mu=\{\fop,\sop\}$, we additionally fill the area bounded to the bottom and to the left by these two half-lines with the same color.

Then, exactly those monomials $\fop^l \sop^j\kappa_i\in \Mon(\cF)$ lie in the $\mu$-cone of $b$, for which the corresponding point $(l,j)$ in the $\kappa_i$-sheet is colored.

\end{rem}

\subsection{Janet basis}

\begin{defn}\label{def:janet-basis}
A \emph{Janet basis} for a submodule $\Dgen{g_1,\ldots,g_s}\subseteq \cF$ is a finite set $T=\{ (b_1,\mu_1),\ldots,(b_r,\mu_r)\}\subseteq \cF\times \cP(\{\fop,\sop\})$ satisfying the following conditions.
\begin{enumerate}
\item $\Dgen{b_1,\ldots,b_r}=\Dgen{g_1,\ldots,g_s}$,
\item \[  \left\{ \lm(f) \mid f\in \Dgen{b_1,\ldots, b_r} \right\} = \biguplus_{i=1}^r \Mon(\mu_i)\lm(b_i)\]
where the operator symbol $\biguplus$ means that the union of sets is a disjoint union.
\end{enumerate}
\end{defn}

\begin{rem}
The definition of a Janet basis is given differently in \cite{dr:faepde}, namely by means of the result of Algorithm \ref{algo:janet-basis}. However, the properties given in the definition here, are the properties needed for the following theorem. 
When we will need the stronger form of the Janet basis obtained from the algorithm, we will explicitly state it.
\end{rem}

\begin{thm}[{\cite[Theorem 2.1.43(d)+(c)]{dr:faepde}}]\label{thm:robertz-K-basis}
Let $T=\{ (b_1,\mu_1),\ldots,(b_r,\mu_r)\}$ be a Janet basis for $\Dgen{g_1,\ldots,g_s}\subseteq \cF$.
\begin{enumerate}
\item \label{item:same-classes} Two elements $f_1,f_2\in \cF$ represent the same element in $\cF/\Dgen{g_1,\ldots,g_s}$, if and only if $\NF(f_1,T)=\NF(f_2,T)$.
In particular, an element $f\in \cF$ is contained in $\Dgen{g_1,\ldots,g_s}$ if and only if $\NF(f,T)=0$.
\item \label{item:K-basis} The residue classes of the elements of
\[ S:= \Mon(\cF) \setminus \bigcup_{i=1}^r \Mon(\mu_i)\lm(b_i).\]
 are a $K$-basis for the factor module $\cF/\Dgen{g_1,\ldots, g_s}$.
\end{enumerate}
\end{thm}

An immediate consequence of this theorem is the following corollary.

\begin{cor}\label{cor:S-basis-of-normal-forms}
Let $T=\{ (b_1,\mu_1),\ldots,(b_r,\mu_r)\}$ be a Janet basis for $\Dgen{g_1,\ldots,g_s}\subseteq \cF$. Let
\[ S:= \Mon(\cF) \setminus \bigcup_{i=1}^r \Mon(\mu_i)\lm(b_i),\quad
\text{and let}
\quad  \rN:= \left\{ f\in \cF \mid f=\NF(f,T) \right\}\]
be the $K$-subspace of $\cF$ of normal forms. Denote by $\pr:\cF\to  \cF/\Dgen{g_1,\ldots, g_s}$ the canonical projection. Then
\begin{enumerate}
\item $S$ is a $K$-basis of $\rN$.
\item The restriction $\pr|_{\rN}:\rN\to \cF/\Dgen{g_1,\ldots, g_s}$ is a $K$-linear isomorphism.
\end{enumerate}
\end{cor}

\subsection{Janet algorithm}

For constructing a Janet basis, we need the notion of \emph{Janet decomposition}.
The general definition is given in \cite[Algorithm 2.1.6]{dr:faepde}. In our setting, this definition simplifies to the following.

\begin{defn}\label{def:janet-decomposition}
Let $G\subseteq \cF\setminus \{0\}$ be a finite set such that for all $g_1,g_2\in G$, $g_1\ne g_2$, the leading monomial $\lm(g_1)$ is not divisible by $\lm(g_2)$, i.e., $\lm(g_1)\not\in \Mon(\{\fop,\sop\})\lm(g_2)$.
The \emph{Janet decomposition} of $G$ is the finite set $J\subseteq \cF\times \cP(\{\fop,\sop\})$ obtained by the following process:
\begin{enumerate}
\item Order the elements in $G$ to obtain a numeration $\{g_1,\ldots, g_s\}$ descending by their leading monomials (i.e.,~$\lm(g_{j+1})\prec \lm(g_j)$).
\item For each $j=1,\ldots, s$ do the following:\\
If for some $i\in \{1,\ldots, d\}$, $g_j$ is the first element in $G$ with $\lm(g_j)\in \Mon(\{\fop,\sop\})\kappa_i$, add the pair $(g_j, \{\fop,\sop\})$ to $J$.
If this is not the case, add the pairs $(\fop^k g_j, \{\sop\})$ to $J$ for all $0\leq k< \deg_\fop(g_{j-1})-\deg_\fop(g_{j})$. 
\end{enumerate}
\end{defn}

\begin{rem}
The Janet decomposition $J$ of $G$ is a special cone decomposition, that means it satisfies
\[  \biguplus_{(b,\mu)\in J} \Mon(\mu)\lm(b) = \bigcup_{g\in G} \Mon(\{\fop,\sop\})\lm(g) \]
\end{rem}

The following example illustrates the idea.

\begin{exmp}\label{exmp:running-example}
Let $\theta\in K$ be an element such that $\fop\cdot \theta=\theta^q\cdot \fop$ and $\sop\cdot \theta=\theta\cdot \sop$.
Consider $\cF=\cD\kappa_1\oplus \cD\kappa_2$, and
\begin{align*}
  g_1 &=\fop^2\kappa_1+(-\sop^2+\theta)\kappa_2,\qquad g_2=(-\sop^2+\theta)\kappa_1+\kappa_2,\\ 
  g_3 &=\left(\fop^2-(\sop^2-\theta)(\sop^2-\theta^{q^2})\right)\kappa_2.
\end{align*}
The generated cones are depicted in the following figure.

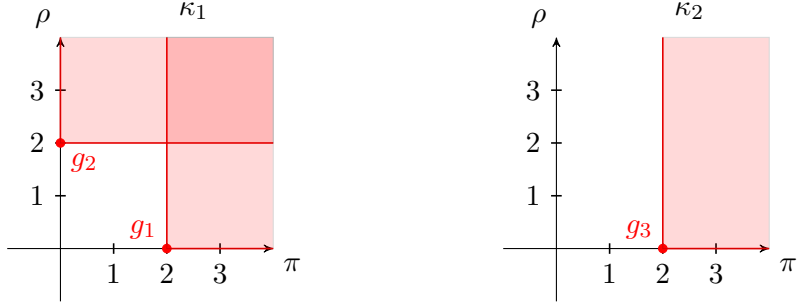
\twosheets{
	\grid{\fop}{\sop}{1}
	\point{g_2}{below right}{0}{2}
	\point{g_1}{above left}{2}{0}
    \cone{2}{0} \cone{0}{2}
}{
	\grid{\fop}{\sop}{2}
	\point{g_3}{above left}{2}{0}
    \cone{2}{0}
}

We have already ordered the elements by descending leading monomials, since $\fop^2\kappa_1\succ \sop^2\kappa_1\succ \fop^2\kappa_2$.
For obtaining the Janet decomposition, we start with $g_1$, and add the pair $(g_1,\{\fop,\sop\})$ to $J$. Next, the leading monomial of $g_2$ is also a multiple of $\kappa_1$,
and $\deg_\fop(g_1)-\deg_\fop(g_2)=2-0=2$. Hence, we add the pairs $(g_2, \{\sop\})$ and $(\fop g_2, \{\sop\})$ to $J$. Finally, $g_3$ is the first element for which $\lm(g_3)$ is a multiple of $\kappa_2$, hence we add $(g_3,\{\fop,\sop\})$ to $J$.

Therefore 
\[ J= \bigl\{(g_1,\{\fop,\sop\}), (\fop g_2,\{\sop\}), (g_2,\{\sop\}), (g_3,\{\fop,\sop\}) \bigr\} \]
is the Janet decomposition of $\{g_1,g_2,g_3\}$.

The cones for $J$ look as follows:

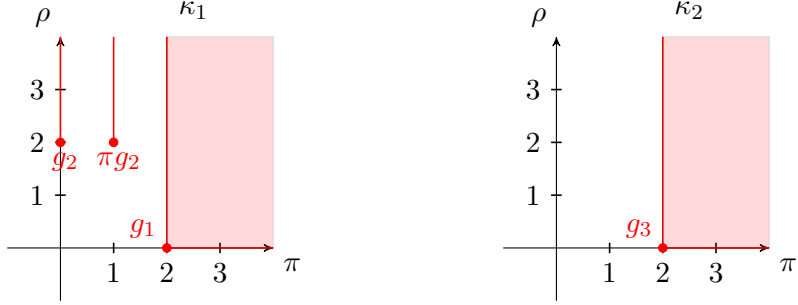
\twosheets{
	\grid{\fop}{\sop}{1}
	\point{\,\,g_2}{below}{0}{2}
	\point{\,\,\fop g_2}{below}{1}{2}
	\point{g_1}{above left}{2}{0}
    \linet{0}{2}
    \linet{1}{2}
    \cone{2}{0} 
}{
	\grid{\fop}{\sop}{2}
	\point{g_3}{above left}{2}{0}
    \cone{2}{0}
}

\end{exmp}

We can now state an algorithm to compute a Janet basis which is a simplified version of \cite[Algorithm 2.1.42]{dr:faepde}, and we refer to \cite[Theorem 2.1.43]{dr:faepde} for a proof that the result is indeed a Janet basis.

\begin{algo}[Janet Algorithm]\label{algo:janet-basis} \
\begin{description}
\item[Input] A finite set $\{g_1,\ldots, g_s\}\subseteq \cF\setminus \{0\}$.
\item[Output] A Janet basis $T=\{ (b_1,\mu_1),\ldots,(b_r,\mu_r)\}$ for $\Dgen{g_1,\ldots,g_s}$.
\item[Algorithm] Set $G:=\{g_1,\ldots, g_s\}$, and repeat the following (1)--(4):
\begin{enumerate}
\item auto-reduction: As long as there is $g\in G$ such that $h:=\NF(g, G\setminus\{g\})\ne g$, replace $g$ by $h$.
\item compute the Janet decomposition $J$ of $G$,
\item set $P:= \bigl\{ \NF(\nu\cdot b,J) \bigmid (b,\mu)\in J, \nu\in \{\fop,\sop\}\setminus \mu \bigr\}\setminus \{0\}$,
\item replace $G$ by $\{ b \mid (b,\mu)\in J\} \cup P$,
\end{enumerate}
until $P=\emptyset$.\\
Return $J$.
\end{description}
\end{algo}

\begin{exmp}\label{exmp:running-example-2}
Let's continue with the example above, $G=\{g_1,g_2,g_3\}$, and compute a Janet basis for it.

Since the leading monomials of the elements in $G$ don't divide any other monomial, the normal forms in step (1) are all equal to the elements themselves. Hence nothing changes there. As stated in the example, the computed Janet decomposition of $G$ is 
\[ J= \bigl\{(g_1,\{\fop,\sop\}), (\fop g_2,\{\sop\}), (g_2,\{\sop\}), (g_3,\{\fop,\sop\}) \bigr\} \]
so in step (3), we have to consider only the normal forms $\NF(\fop\cdot \fop g_2,J)$ and $\NF(\fop \cdot g_2, J)$. The second one is trivially $0$, since $(\fop g_2,\{\sop\})\in J$. For the first one, we compute:
\[ \fop^2 g_2= -\sop^2g_1+\theta^{q^2}g_1+g_3. \]
Hence, also $\NF(\fop^2 g_2, J)=0$. This means that $P=\emptyset$, and $J$ is indeed a Janet basis for $G$.
\end{exmp}

\section{Finite generation and relations}\label{sec:finite-generation-and-relations}

In this section, we prove our main theorems which will be applied to the questions concerning Anderson $t$-modules, their $t$-motives, and $t$-comotives.

In addition to the notation in the previous sections, we assume that $M$ is finitely generated free as a $K\sp{\fop}$-module, and $(\bkappa_1,\ldots, \bkappa_d)$ is a $K\sp{\fop}$-basis of $M$.
So if the $\sop$-action is given in this basis by a matrix $D\in \Mat_{d\times d}(K\sp{\fop})$, i.e.,~in matrix notation
\begin{equation} \label{eq:sop-action} 
\sop\cdot \svect{\bkappa}{d} = D\cdot \svect{\bkappa}{d},
\end{equation}
the obvious generators in the presentation \eqref{eq:presentation-of-M}, are $\{p_1,\ldots, p_d\}$ where for $i=1,\ldots, d$:
\[  p_i := \sop\kappa_i - \sum_{j=1}^d D_{ij}\kappa_j\in \cF. \]

As in the previous section, we use the lexicographical monomial order on 
$\cD=\Kpolys$ with $\sop\prec \fop$, and extend this order to $\cF$ via the position-over-term order.
Further, we let $J=\{(b_1,\mu_1),\ldots, (b_s,\mu_s)\}$ be a Janet basis for $\Dgen{p_1,\ldots, p_d}$ with respect to this order.

We also let $\pr:\cF\to M$ be the canonical projection, i.e.,~the $\Kpolys$-linear map sending $\kappa_i$ to $\bkappa_i$, and let
\[  \rN:= \left\{ f\in \cF \mid f=\NF(f,J) \right\}\]
be the $K$-subspace of normal forms.
We recall from Corollary \ref{cor:S-basis-of-normal-forms} that the set 
 \[ S:= \Mon(\cF) \setminus \bigcup_{i=1}^r \Mon(\mu_i)\lm(b_i)\]
is a $K$-basis of $\rN$, and that the restriction $\pr|_{\rN}:\rN\to M$ is a $K$-linear isomorphism. 
In formulas, we sometimes emphasize this fact be writing $\bar{f}$ instead of $\pr(f)$ for a normal form $f\in \rN$.

\begin{defn}\label{def:quantities}
Depending on $J$, we define the following quantities (with the convention that the infimum of an empty set is $\infty$).
For each $i=1,\ldots, d$:
\begin{align*}
\sfn_i&:=\inf \{  n\in \NN \mid \exists (b,\mu)\in J \text{ s.t. }  \lm(b)=\fop^n\kappa_i \}\ \in \NN\cup \{\infty \}, \\
\sfm_i&:=\inf \{  n\in \NN \mid \exists (b,\mu)\in J,l\geq 0\text{ s.t. }  \lm(b)=\fop^n \sop^l\kappa_i \} \in \NN\cup \{\infty \}, \\
\text{as well as}& \\
\Wgen &:= \{ \fop^{j}\kappa_i \mid i=1,\ldots, d, 0\leq j<\sfn_i \} \ \subseteq \cF, \\
\Wind &:= \{ \fop^{j}\kappa_i \mid i=1,\ldots, d, 0\leq j<\sfm_i \} \ \subseteq \cF.
\end{align*}
\end{defn}

Obviously, we have $\sfm_i\leq \sfn_i$ for all $i=1,\ldots, d$, and hence also $
\Wind \subseteq \Wgen$.

\begin{prop}\label{prop:properties-of-quantities} \
\begin{enumerate}
\item \label{item:restriction-of-pi} The restriction of $\pr$ to $\Wgen$ is injective.
\item \label{item:wgen} $\bWgen:=\pr(\Wgen)$ is generating $M$ as a $K\sp{\sop}$-module.
\item \label{item:m_i-finite} For all $i=1,\ldots, d$, we have $\sfm_i<\infty$.
\item \label{item:wind} $\bWind:=\pr(\Wind)$ is a maximally $K\sp{\sop}$-linearly independent subset of $M$.
\end{enumerate}
\end{prop}

\begin{proof} 
\noindent \eqref{item:restriction-of-pi}: By definition of $\sfn_i$, we have $\Wgen\subseteq S$, and by Corollary \ref{cor:S-basis-of-normal-forms}, the restriction of $\pr$ to $S$ is injective.

\noindent \eqref{item:wgen}: By definition of the $\sfn_i$ and of $\Wgen$, we have
\[  S \subseteq \{ \fop^j \sop^l\kappa_i\mid i=1,\ldots, d,\, 0\leq j< \sfn_i,l\geq 0 \}=\Mon(\{\sop\})\cdot \Wgen. \]
Since $\pr(S)$ is a $K$-basis of $M$, $M$ is generated as a $K\sp{\sop}$-module by $\bWgen$.

\noindent \eqref{item:m_i-finite}: Suppose to the contrary that some $\sfm_i=\infty$. Then 
\[ S\supseteq \{ \fop^k \sop^l\kappa_i \mid k,l\geq 0 \}.\]
Hence, the elements 
$\fop^k \sop^l\bkappa_i$, $k,l\geq 0$, are $K$-linearly independent in $M$, and therefore, the elements $\sop^l\bkappa_i$, $l\geq 0$, are $K\sp{\fop}$-linearly independent. This, however, is not possible, since by assumption, $M$ is finitely generated over the left principal ideal domain $K\sp{\fop}$.

\noindent \eqref{item:wind}: By definition of $\sfm_i$ and $\Wind$, we have
\[  S \supseteq \{ \fop^j \sop^l\kappa_i\mid i=1,\ldots, d,\, 0\leq j< \sfm_i,l\geq 0 \}=\Mon(\{\sop\})\cdot \Wind. \]
Since, $\pr(S)$ is a $K$-basis of $M$, the set $\Mon(\{\sop\})\cdot \bWind$ is $K$-linearly independent, and hence $\bWind$ is a $K\sp{\sop}$-linearly independent set.

For showing maximality, suppose to the contrary that there is $g\in \cF$ such that $\pr(g)$ is $K\sp{\sop}$-linearly independent from $\bWind$. Of all those elements, we choose one with $\lm(g)$ being minimal. Now let $i,j,l$ be such that $\lm(g)=\fop^j\sop^l\kappa_i$.

If we had $j< \sfm_i$, then we could subtract a $K\sp{\sop}$-multiple of some element of $\Wind$, and obtain another such $g$ with smaller leading monomial. Hence, $j\geq \sfm_i$.

By definition of $\sfm_i$, there is $k\geq 0$, and $(b,\mu)\in J$ such that
$\fop^j\sop^{l+k}\kappa_i$ lies in the $\{\fop,\sop\}$-cone generated by $b$, and even more, we can choose $(b,\mu)\in J$ such that
\[ \sop^{k}\lm(g)=\fop^j\sop^{l+k}\kappa_i\in \Mon(\mu)\lm(b).\]

Hence, the leading monomial of $\NF(\sop^{k}g, J)$ is strictly smaller than $\sop^{k}\lm(g)$ which means that $f=\sop^{k}g-\NF(\sop^{k}g, J)\in \cF$ is a non-zero element with $\lm(f)=\sop^{k}\lm(g)$. Further by Theorem \ref{thm:robertz-K-basis}\eqref{item:same-classes}, we have $\pr(f)=0$. 

Due to our chosen monomial order, we can write
\[ f = f_1(\sop)g + r, \]
with a non-trivial polynomial $f_1(\sop)\in K\sp{\sop}$\footnote{Obviously, $f_1(\sop)$ is even monic and of degree $k$, but that's not important for the proof.}, and $r\in \cF$ with $\lm(r)<\lm(g)$.

Since $\pr(g)$ was $K\sp{\sop}$-linearly independent from $\bWind$, also $\pr(f_1(\sop)g)=f_1(\sop)\pr(g)$ is $K\sp{\sop}$-linearly independent from $\bWind$, and so is
\[ \pr(r)=\pr(f-f_1(\sop)g)=-f_1(\sop)\pr(g). \]
This, however, contradicts our assumption that $g$ was such an element with minimal leading monomial.
\end{proof}

We can now prove our first main theorem.


\begin{thm}\label{thm:criterion-for-finiteness}
With the definition and notation from above.
\begin{enumerate}
\item \label{item:M-finite} $M$ is finitely generated over $K\sp{\sop}$ $\Leftrightarrow$ $\forall i=1,\ldots, d$, we have $\sfn_i<\infty$.
\item \label{item:rational-dimension}
Let $K(\sop)$ denote the skew field of fractions of $K\sp{\sop}$. The space $K(\sop)\otimes_{K\sp{\sop}} M$ is a finite dimensional $K(\sop)$-vector space of dimension
\[ \dim_{K(\sop)}(K(\sop)\otimes_{K\sp{\sop}} M) = \sum_{i=1}^d \sfm_i. \]
\item \label{item:sop-rank} If $M$ is $K\sp{\sop}$-finitely generated, and $\gamma_\sop:K\to K$ is bijective, then $M$ is a free $K\sp{\sop}$-module of rank 
 \[ \rk_{K\sp{\sop}}(M) = \sum_{i=1}^d \sfm_i.\]
\end{enumerate}
\end{thm}

\begin{proof}
\begin{enumerate}
\item[\eqref{item:M-finite}:] If $\sfn_i<\infty$ for all $i=1,\ldots, d$, then $\Wgen$ is finite. Hence by Proposition~\ref{prop:properties-of-quantities}\eqref{item:wgen}, $M$ is finitely generated as $K\sp{\sop}$-module. 

On the other hand, assume that for some $i=1,\ldots, d$, we have $\sfn_i=\infty$. 
By definition of $\sfn_i$, this implies that for all $n\geq 0$, $\fop^n\kappa_i\in S$, and hence $\NF(\fop^n\kappa_i, J)=\fop^n\kappa_i$.

Further, due to our chosen monomial order, $\fop^j\sop^l \kappa_i \prec \fop^n\kappa_i$ for all $j<n$ and $l\geq 0$. In particular, for
any element $f=\sum_{j=0}^{n-1} \sum_{l=0}^{L_j} x_{jl} \fop^j\sop^l \kappa_i$ in the subspace $\Ksopgen{\fop^{j}\kappa_i \mid j<n}\subset \cF$, we have
\[  \lm\left( \NF(f,J)\right) \preceq \lm(f)\prec \fop^n\kappa_i = \NF(\fop^n\kappa_i, J). \]
which implies $\NF(f,J) \ne \NF(\fop^n\kappa_i, J)$. Hence, by Theorem \ref{thm:robertz-K-basis}\eqref{item:same-classes}, $\pr(f)\ne \pr(\fop^n\kappa_i)=\fop^n\bkappa_i$.

This show that for all $n\geq 0$, $\fop^{n}\bkappa_i\notin \Ksopgen{\fop^{j}\bkappa_i \mid j<n}\subseteq M$.
Hence, the $K\sp{\sop}$-submodules $\Ksopgen{\fop^{j}\bkappa_i \mid j<n}$, $n\geq 0$, form an infinitely ascending chain of submodules of $M$ which is impossible, if $M$ is finitely generated, since $K\sp{\sop}$ is Noetherian. Hence, $M$ is not finitely generated as $K\sp{\sop}$-module.
\item[\eqref{item:rational-dimension}:] By Proposition \ref{prop:properties-of-quantities}\eqref{item:wind}, the set $\bWind=\pr(\Wind)$ is a maximally $K\sp{\sop}$-linearly independent subset of $M$. Therefore, it is a maximally $K(\sop)$-linearly independent subset of $K(\sop)\otimes_{K\sp{\sop}} M$, hence a basis.
By Proposition \ref{prop:properties-of-quantities}\eqref{item:m_i-finite}, the set $\Wind$ is finite, and its cardinality is $\sum_{i=1}^d \sfm_i$, which proves the claim.
\item[\eqref{item:sop-rank}:] If $M$ is finitely generated, and $\gamma_\rho$ is an isomorphism, then by Proposition \ref{prop:observations-on-M}\eqref{item:freeness}, it is not only finitely generated, but even a free $K\sp{\sop}$-module. Hence, 
\[ \rk_{K\sp{\sop}}(M)=\dim_{K(\sop)} \left( K(\sop)\otimes_{K\sp{\sop}} M \right) = \sum_{i=1}^d \sfm_i \]
by Part \eqref{item:rational-dimension}.\qedhere
\end{enumerate}
\end{proof}

We are next going to describe the $\fop$-action on the $K\sp{\sop}$-generators $\bWgen$ of $M$, as well as the $K\sp{\sop}$-relations among the elements of $\bWgen$.
For obtaining this nice description, we need a Janet basis that we obtain from Algorithm \ref{algo:janet-basis}.

\begin{lem}\label{lem:description-of-J}
Assume that $\sfn_i<\infty$ for all $i=1,\ldots, d$, and that $J$ is obtained as the result of Algorithm \ref{algo:janet-basis}.
For all $i=1,\ldots, d$, let $(b_i,\mu_i)\in J$ be the pair such that $\lm(b_i)=\fop^{\sfn_i}\kappa_i$, and set
\[  \Jtop:=\{ (b_i,\mu_i) \mid i=1,\ldots, d\},\quad \text{and}\quad \Jlow:=J\setminus \Jtop. \]
Then
\begin{enumerate}
\item \label{item:mu_i} for all $i=1,\ldots, d$, we have $\mu_i=\{ \fop,\sop\}$.
\item \label{item:pairs_in_J_1} For all $i=1,\ldots, d$, for all $\sfm_i\leq j<\sfn_i$, there exist unique $l>0$ and $(b,\mu)\in \Jlow$ such that $\lm(b)=\fop^j\sop^l\kappa_i$. For those $(b,\mu)$, one has $\mu=\{\sop\}$.
\item \label{item:description_of_J_1} $\Jlow$ consists exactly of the pairs $(b,\mu)$ given in \eqref{item:pairs_in_J_1}.
\end{enumerate}
\end{lem}

\begin{proof}
Assume that $G$ is the auto-reduced set obtained after step (1) in the last run of the loop in Algorithm \ref{algo:janet-basis}, so that $J$ equals the Janet decomposition of $G$. We start by showing \eqref{item:mu_i}. Let $i\in \{1,\ldots, d\}$.

If there wasn't any $b\in G$ with $\lm(b)\in \Mon(\{\fop\})\kappa_i$, then also $J$ wouldn't contain such an element, contradicting our assumption $\sfn_i<\infty$. So let $b\in G$ with $\lm(b)=\fop^n\kappa_i$ be the one with least $n$. Since $G$ is auto-reduced, all other elements $\tilde{b}\in G$ with $\lm(\tilde{b})\in \Mon(\{\fop,\sop\})\kappa_i$ have $\deg_{\fop}(\lm(\tilde{b}))<n$. In particular, $b$ is the only element in $G$ whose leading monomial lies in $\Mon(\{\fop\})\kappa_i$. By construction of the Janet decomposition, it remains the only such element, and the corresponding $\mu$ equals $\{\fop,\sop\}$. So $b=b_i$ and $\mu_i=\{\fop,\sop\}$, showing \eqref{item:mu_i}.

Part \eqref{item:pairs_in_J_1} and Part \eqref{item:description_of_J_1} are then also directly obtained by the construction of the Janet decomposition of $G$, again taking into account that $G$ is auto-reduced.
\end{proof}

\begin{rem}\label{rem:monomials-of-b-i}
As the Janet basis in the previous lemma is obtained as the Janet decomposition of the auto-reduced set $G$, and all $b_i$ were already elements of $G$, we also obtain that for all $i=1,\ldots, d$, apart from the leading monomial $\lm(b_i)$, all other monomials that occur in $b_i$ are in $S$. 

The same holds for $(b,\mu)\in \Jlow$, if $b\in G$. However, elements that were added in the process of the Janet decomposition, i.e.,~those of the form $\fop^k g$ with $g\in G$ might not have this property.
For those $(b,\mu)$, we instead use $\tilde{b}=\NF(b,J\setminus \{(b,\mu)\})$ in the next proposition. Then apart from the leading monomial $\lm(\tilde{b})$ which equals $\lm(b)$, all other monomials that occur in $\tilde{b}$ are in $S$. 
\end{rem}

\begin{thm}\label{thm:relations-and-fop-action}
Assume that $\sfn_i<\infty$ for all $i=1,\ldots, d$, and that $J$ is obtained as the result of Algorithm \ref{algo:janet-basis}.
Let \footnote{See Remark \ref{rem:monomials-of-b-i} for the choice of $\Blow$.}
\begin{align*}
 \Btop &:=\{ b\mid (b,\mu)\in \Jtop\}=\{b_1,\ldots,b_d\}, \quad \text{and}\\ 
 \Blow &:=\bigl\{ \NF(b, J\setminus \{(b,\mu)\}) \bigmid (b,\mu)\in \Jlow\bigr\} 
\end{align*}
with $\Jtop$, $\Jlow$, and $b_1,\ldots,b_d$ as in Lemma \ref{lem:description-of-J}. Then
\begin{enumerate}
\item for $w\in \Wgen$, i.e.,~$w=\fop^j\kappa_i$ for some $i=1,\ldots, d$, and $0\leq j< \sfn_i$ one has
\[   \NF(\fop w,J)= \partdef{ \fop w & \text{if } j<\sfn_i-1\\ \fop w-\frac{1}{\lc(b_i)}b_i & \text{if } j=\sfn_i-1\,.  } \]
\item A generating set for the $K\sp{\sop}$-relations among the elements in $\bWgen$ is given by writing the elements in $\Blow$ as $K\sp{\sop}$-linear combinations of the elements in $\Wgen$.
\end{enumerate}
\end{thm}

\begin{proof}
\begin{enumerate}
\item Let $w=\fop^j\kappa_i\in \Wgen$. If $j<\sfn_i-1$, then $j+1<\sfn_i$, and $\fop w=\fop^{j+1}\kappa_i\in S$. Hence it equals its normal form. If $j=\sfn_i-1$, then $\fop^{j+1}\kappa_i=\fop^{\sfn_i}\kappa_i=\lm(b_i)$. Hence, the first step in the reduction process is to subtract $\frac{1}{\lc(b_i)}b_i$. Further from Remark \ref{rem:monomials-of-b-i}, we see that all monomials of $\fop w-\frac{1}{\lc(b_i)}b_i$ are elements of $S$. Therefore,  $\fop w-\frac{1}{\lc(b_i)}b_i$ is the normal form of $\fop w$.

\item By the main property of the Janet basis, all $K\sp{\sop}$-relations among elements in $K\sp{\sop}\bWgen$ are given by the set $\{ f-\NF(f,J) \mid f\in K\sp{\sop}\Wgen \}$. By the definition of the normal forms and of $\Btop$ and $\Blow$, one has
\begin{equation}\label{eq:f-minus-NF-f}
  f- \NF(f,J)= \sum_{(b,\mu)\in J} \tilde{c}_b b = \sum_{b\in \Btop\cup \Blow} c_b b 
\end{equation}
for appropriate $\tilde{c}_b,c_b\in \Kpolys$, where $\tilde{c}_b\in K\sp{\sop}$ if $(b,\mu)\in \Jlow$ and $c_b\in  K\sp{\sop}$ if $b\in \Blow$.
The latter form is obtained from the former by replacing the $b$'s from $\Jlow$ by the linear combination obtained when computing their normal forms.

Since we only consider $f\in K\sp{\sop}\Wgen$, each monomial $\fop^j\sop^l\kappa_i$ occurring in $f$ satisfies $j<\sfn_i$. The same is true for all normal forms, hence for $\NF(f,J)$.
For given $i=1,\ldots, d$, the only monomials $\fop^j\sop^l\kappa_i$ in $\sum_{b\in \Btop\cup \Blow} c_b b$ which satisfy $j\geq \sfn_i$ stem from 
the term $c_{b_i}b_i$, and hence $c_{b_i}=0$ for all $i=1,\ldots, d$, i.e.,~$c_b=0$ for all $b\in \Btop$.

This shows that all differences $f-\NF(f,J)$ are $K\sp{\sop}$-linear combinations of the elements in $\Blow$, proving the claim.\qedhere
\end{enumerate}
\end{proof}

\begin{rem}\label{rem:description-if-sfn_i=infty}
In the previous theorem, we focused on the case that $\sfn_i<\infty$ for all $i=1,\ldots,d$, i.e., that $M$ is a finitely generated $K\sp{\sop}$-module, because this is the interesting case for our applications. However, a variation of Lemma \ref{lem:description-of-J} and Theorem \ref{thm:relations-and-fop-action} will also lead to an explicit description of the relations and the $\fop$-action on $\bWgen$.
We will give an example thereof in Example \ref{exmp:quasi-periodic-extension}.
\end{rem}

\begin{rem}\label{rem:finding-a-basis}
If we assume that $M$ is finitely generated and $\gamma_\sop$ is an isomorphism, we obtain from Proposition \ref{prop:observations-on-M}\eqref{item:freeness}, that $M$ is a free $K\sp{\sop}$-module. From the finite generating set  $\bWgen$ and the finite relations obtained from all $b\in \Blow$, one can obtain a basis for $M$. 
In order to achieve this, one applies the elementary divisor algorithm to this setting, which also works in this non-commutative case, as $K\sp{\sop}$ is both left and right principal (cf.~\cite[Ch.~3.1 \& 3.7]{nj:tr}).\footnote{This is exactly the point where one needs that $\gamma_\sop$ is an isomorphism. Otherwise $K\sp{\sop}$ would not be right principal.}

For the convenience of the reader, we give more details. Let $\Wgen=\{w_1,\ldots, w_s\}$, 
write all $b\in \Blow$ as $K\sp{\sop}$-linear combinations of the elements in $\Wgen$, and build the matrix $B\in \Mat_{r\times s}(K\sp{\sop})$ whose rows are the coefficient vectors of these linear combinations. Here we denote $r=\# \Blow$.
Diagonalize the matrix via row and column operations with a bookkeeping of the column operations. The diagonalizing amounts to finding matrices $U\in \GL_{r}(K\sp{\sop})$ and $V\in \GL_{s}(K\sp{\sop})$ such that (in block matrix form)
\[  UBV = \begin{pmatrix}
 H & 0 \\ 0 & 0 
\end{pmatrix}, \]
where $H$ is a square diagonal matrix of some size $s_0\leq \min\{ r,s\}$ with non-zero entries on the diagonal, and the $0$'s denote zero matrices of the appropriate sizes.

As $M$ is torsion-free, it turns out that the diagonal entries of $H$ are units, i.e., we can achieve that $H=\one_{s_0}$ is the identity matrix.

Then the set $\{e_1,\ldots, e_{s-s_0} \}$ where 
\[ e_{i}:= \sum_{j=1}^s (V^{-1})_{s_0+i,j} \bar{w}_j \]
forms a $K\sp{\sop}$-basis of $M$.

\medskip

Using the transformation matrix $V$, one also easily computes the $\fop$-action on $\{e_1,\ldots, e_{s-s_0} \}$ from that on $\bWgen$. Namely, let in matrix notation
\[  \fop \svect{\bar{w}}{s} = C \svect{\bar{w}}{s} \]
with $C\in \Mat_{s\times s}(K\sp{\sop})$. Then with $s_1=s-s_0$,
\[  \fop \svect{e}{s_1} = \tilde{C} \svect{e}{s_1} \]
where $\tilde{C}$ is the lower right $(s_1\times s_1)$-submatrix of $\gamma_\fop(V^{-1})CV$.\footnote{Here, $\gamma_\fop(V^{-1})$ means the entrywise application of $\gamma_\fop:K\sp{\sop}\to K\sp{\sop}$ (the $\sop$-linear extension of the original $\gamma_\fop$).}

We leave the verification of this computation as an exercise to the reader.
\end{rem}

\begin{rem}\label{rem:changing-coordinate-system}
The key point in the results of this section is the chosen monomial order in $\Kpolys$, and that we use a position-over-term order for $\cF$. However, we are free to change the basis $\{\kappa_1,\ldots,\kappa_d\}$, as it suits best. 
As such a change does also change the order on $\cF$, if can affect the run-time of the algorithm in special cases (see Example \ref{exmp:quasi-periodic-extension}).
\end{rem}

We end this section by illustrating it with our running Example \ref{exmp:running-example}.

\begin{exmp}\label{exmp:running-example-3}
In Example \ref{exmp:running-example-2}, we computed the Janet basis
\[ J= \bigl\{(g_1,\{\fop,\sop\}), (\fop g_2,\{\sop\}), (g_2,\{\sop\}), (g_3,\{\fop,\sop\}) \bigr\}, \]
where
\begin{align*}  
g_1&=\fop^2\kappa_1+(-\sop^2+\theta)\kappa_2,&\quad g_2&=(-\sop^2+\theta)\kappa_1+\kappa_2,\\
\fop g_2&=(-\fop\sop^2+\theta^q\fop)\kappa_1+\fop \kappa_2,& \quad g_3&=\left(\fop^2-(\sop^2-\theta)(\sop^2-\theta^{q^2})\right)\kappa_2. 
\end{align*}

We readily see that in this example, we have
\begin{align*}
\sfn_1 &= 2,\, \sfm_1=0,\quad \sfn_2 =\sfm_2=2,\\
\Wgen &= \{ \kappa_1,\fop\kappa_1, \kappa_2,\fop\kappa_2\}, \\
\Wind &= \{ \kappa_2,\fop\kappa_2 \}.
\end{align*}
So by Theorem \ref{thm:criterion-for-finiteness}, and Proposition \ref{prop:properties-of-quantities}, $M$ is a finitely generated $K\sp{\sop}$-module, and
$\bWgen = \{ \bkappa_1,\fop\bkappa_1, \bkappa_2,\fop\bkappa_2\}$ is a generating set for $M$.

Further, we have
\[ \Btop = \{  g_1, g_3 \},\quad \text{and}\quad \Blow=\{ \fop g_2, g_2\}.\]

The $K\sp{\sop}$-relations are hence given by
\begin{align*}
0 &= \pr(\fop g_2) = (-\sop^2+\theta^q)\fop \bkappa_1 +\fop\bkappa_2,\quad \text{and}\\ 
0 &= \pr(g_2) = (-\sop^2+\theta^q) \bkappa_1 +\bkappa_2. 
\end{align*}
So $\fop\bkappa_2$ and $\bkappa_2$ are $K\sp{\sop}$-multiples of $\fop \bkappa_1$ and $\bkappa_1$, respectively. Therefore, $\{ \bkappa_1,\fop\bkappa_1\}$ is a $K\sp{\sop}$-basis for $M$, in accordance with Part \eqref{item:sop-rank} of Theorem \ref{thm:criterion-for-finiteness}, stating that the rank is $\sfm_1+\sfm_2=2$.

We finally compute the $\fop$-action with respect to the basis $\{ e_1:=\bkappa_1,e_2:=\fop\bkappa_1\}$ using Theorem \ref{thm:relations-and-fop-action}:
\begin{align*}
\fop e_1 &= \fop\bkappa_1 = e_2,\\
\fop e_2 &= \pr(\fop^2\bkappa_1 - g_1)=(\sop^2-\theta)\bkappa_2
= (\sop^2-\theta)(\sop^2-\theta^q) \bkappa_1\\
&= (\sop^2-\theta)(\sop^2-\theta^q)e_1.
\end{align*}
\end{exmp}

\section{Anderson modules, motives, and comotives} \label{sec:t-modules}

\subsection{Basic setting}

From now on, let $\Fq$ be the finite field with $q$ elements, and let 
$K$ be a field containing $\Fq$.
Further, let $\Fq[t]$ be a polynomial ring over $\Fq$ in an indeterminate $t$ which is linearly independent to $K$, and $\ell:\Fq[t]\to K$ a homomorphism of $\Fq$-algebras. We set $\theta:=\ell(t)$.

We denote by $K\{\tau\}$ the skew polynomial ring
\[  K\sp{\tau}:=\left\{ \sum_{i=0}^n x_i\tau^i \,\middle|\, n\geq 0, x_i\in K\right\} \]
with multiplication uniquely given by additivity and the rule
\[ \tau \cdot x = x^q\cdot \tau, \]
for all $x\in K$, i.e., 
\[ \left( \sum_{i=0}^n x_i\tau^i\right) \cdot \left(\sum_{j=0}^m y_j\tau^j\right)
= \sum_{k=0}^{n+m} \left( \sum_{i=0}^k x_i\cdot (y_{k-i})^{q^i} \right) \cdot \tau^k.  \]
This ring equals the ring $\End_{\grp,\Fq}(\Ga)$ of $\Fq$-linear group endomorphisms of $\Ga$ by identifying $\tau$ with the $q$-th power Frobenius map, and $x\in K$ with scalar multiplication by~$x$.

If $K$ is perfect, i.e.,~the $q$-th power map is an isomorphism, we denote by $K\sp{\sigma}$ the skew polynomial ring in a variable $\sigma$ subject to the rule 
\[\sigma\cdot x = x^{1/q}\cdot \sigma,\] 
for all $x\in K$.

This ring equals the opposite ring of $K\sp{\tau}$.

The $q$-th power Frobenius map will be extended to $K[t]=K\otimes_\Fq \Fq[t]$ as the \emph{Frobenius twist}, denoted $f\mapsto f^{(1)}$, by applying the $q$-th power to the coefficients in $K$, and leaving $t$ invariant. In the same way, if $K$ is perfect, the inverse Frobenius map will be extended to $K[t]$ as the \emph{inverse Frobenius twist}, denoted $f\mapsto f^{(-1)}$, by applying the $1/q$-th power to the coefficients in $K$, and leaving $t$ invariant.

In the following, we will encounter the tensor product $K\sp{\tau}\otimes_\Fq \Fq[t]$, which is nothing else than the skew polynomial ring $K\sp{\tau,t}$ in the two variables subject to $\tau t=t\tau$, $\tau x = \gamma_\tau(x) \tau$, and $tx=\gamma_t(x) t$ for all $x\in K$, where $\gamma_\tau(x)=x^q$ ($x\in K$) is the Frobenius map, and $\gamma_t=\id_K$ is the identity.

Similarly, if $K$ is perfect, $K\sp{\sigma}\otimes_\Fq \Fq[t]=K\sp{\sigma,t}$ with the automorphisms $\gamma_\sigma(x)=x^{1/q}$ for all $x\in K$ and $\gamma_t=\id_K$.

\subsection{Anderson's objects}

In Anderson's theory, there are three main types of objects. The definitions slightly differ from source to source. We follow \cite{uh-akj:pthshcff} here, and adapt it to our setting.

An \emph{effective Anderson $t$-motive} $(\mot,\tau_\mot)$ over $K$ is a free $K[t]$-module of some finite rank $r$ together with a  homomorphism $\tau_\mot:\mot\to \mot$ that is semilinear with respect to the Frobenius twist\footnote{Frobenius twist semilinearity means that the map $\tau_\mot$ is additive and satisfies $\tau_\mot(f\cdot m)=f^{(1)}\tau_\mot(m)$ for all $f\in K[t]$, and $m\in \mot$.} and such that the ``cokernel'' $\mot/\left(K[t]\cdot \tau_\mot(\mot)\right)$ is annihilated by some power of $t-\theta$.

The existence of the semilinear operator $\tau_\mot$ can equivalently be described by saying that $\mot$ is a left module
over the skew polynomial ring $K\sp{\tau,t}$.

Similarly for $K$ perfect, an \emph{effective Anderson $t$-comotive}\footnote{As mentioned earlier, we avoid the confusing terminology \emph{dual $t$-motive} for this object.} $(\dumot,\sigma_\dumot)$ over $K$ is a free $K[t]$-module of some finite rank $r$ together with a homomorphism $\sigma_\dumot:\dumot\to \dumot$ that is semilinear with respect to the inverse Frobenius twist\footnote{Semilinearity for the inverse Frobenius twist means that the map $\sigma_\dumot$ is additive and satisfies $\sigma_\dumot(f\cdot m)=f^{(-1)}\sigma_\dumot(m)$ for all $f\in K[t]$, and $m\in \dumot$.} and such that the ``cokernel'' $\dumot/\sigma_\dumot(\dumot)$ is annihilated by some power of $t-\theta$.

Parallel to Anderson $t$-motives, the existence of the semilinear operator $\sigma_\dumot$ can equivalently be described by saying that $\dumot$ is a left module over the skew polynomial ring $K\sp{\sigma,t}$.

An \emph{Anderson $t$-module} $(E,\phi)$ over $K$ of dimension $d$ is an $\Fq$-vector space scheme $E$ over $K$ isomorphic to $\Ga^d$ together with a homomorphism of $\Fq$-algebras $\phi:\Fq[t]\to \End_{\grp,\Fq}(E)$ into the ring of $\Fq$-linear group scheme endomorphisms of $E$, such that for all $a\in \Fq[t]$, the endomorphism $\dphi_a$ on $\Lie(E)$ induced by $\phi_a:=\phi(a)$ fulfills the condition that $\dphi_a-\ell(a)$ is nilpotent.

\begin{rem}
The definition of an effective Anderson $t$-motive in \cite[Definition 2.3.1(a) and (c)]{uh-akj:pthshcff} looks different to what we stated here. Hartl and Juschka 
required that the linearization $\tilde{\tau}_\mot:K[t]\otimes_{K[t],\tau} \mot \to \mot, x\otimes m\mapsto x\cdot \tau_\mot(m)$ of $\tau_\mot$ becomes an isomorphism after inverting $t-\theta$. Here the tensor product is via the Frobenius twist on the left hand factor, i.e., $x\otimes ym=xy^{(1)}\otimes m$ for $x,y\in K[t]$, $m\in \mot$.\\
The cokernel of $\tilde{\tau}_\mot$ is $\mot/\left(K[t]\cdot \tau_\mot(\mot)\right)$, i.e., what we called ``cokernel'' above. So the cokernel is annihilated by a power of $t-\theta$ if and only if $\tilde{\tau}_\mot$ becomes surjective after inverting $t-\theta$.
Furthermore, $\tilde{\tau}_\mot$ is a homomorphism of free $K[t]$-modules of the same rank $r$. So if its cokernel is torsion, the image of $\tilde{\tau}_\mot$ also has rank $r$, and this implies that $\tilde{\tau}_\mot$ is injective. Therefore, in our setting, the definition given here is indeed equivalent to the one in \cite{uh-akj:pthshcff}.

A similar argument shows that in our setting the given definition of an effective Anderson $t$-comotive is equivalent to \cite[Definition 2.4.1(a) and (c)]{uh-akj:pthshcff}.
\end{rem}

\subsection{Relation between these objects}\label{subsec:relation-between-Andersons-objects}

For an Anderson $t$-module $(E,\phi)$, the \emph{$t$-motive associated to} $(E,\phi)$ is the $K\sp{\tau,t}=K\sp{\tau}\otimes_\Fq \Fq[t]$-module
\[ \mot(E):=\Hom_{\grp,\Fq}(E,\Ga) \]
of $\Fq$-linear homomorphisms of groups schemes
with left $K\sp{\tau}$-action given by composition with elements in $\End_{\grp,\Fq}(\Ga)=K\sp{\tau}$, and $\Fq[t]$-action given by composition with the elements $\phi_a$ for $a\in \Fq[t]$.

In general, however, $\mot(E)$ might not be finitely generated free as $K[t]$-module. The Anderson $t$-module $(E,\phi)$ and its motive $\mot(E)$ are called \emph{abelian}, if $\mot(E)$ is finitely generated free as $K[t]$-module. In this case, the nilpotence of the endomorphism  $\dphi_a-\ell(a)$ implies that the ``cokernel'' $\mot(E)/\left(K[t]\cdot \tau(\mot(E))\right)$ is annihilated by some power of $t-\theta$. Hence in this case, $\mot(E)$ is indeed an effective Anderson $t$-motive.

In parallel, the \emph{$t$-comotive associated to} $(E,\phi)$ is the right $K\sp{\tau}\otimes_\Fq \Fq[t]$-module
\[ \dumot(E):=\Hom_{\grp,\Fq}(\Ga,E) \]
with right $K\sp{\tau}$-action given by composition with elements in $\End_{\grp,\Fq}(\Ga)=K\sp{\tau}$, and $\Fq[t]$-action given by composition with the elements $\phi_a$ for $a\in \Fq[t]$.
Using the opposite ring $K\sp{\sigma}$ of $K\sp{\tau}$ (in case that $K$ is perfect), $\dumot(E)$ is turned into a left module over $K\sp{\sigma,t}=K\sp{\sigma}\otimes_\Fq \Fq[t]$.
As above,  $\dumot(E)$ might not be finitely generated free as $K[t]$-module in general, and the Anderson $t$-module $(E,\phi)$ and its comotive $\dumot(E)$ are called \emph{coabelian} (or \emph{$t$-finite}), if $\dumot(E)$ is finitely generated free as $K[t]$-module. In this case, the nilpotence of the endomorphism  $\dphi_a-\ell(a)$ implies that $\dumot(E)/\sigma_{\dumot(E)}(\dumot(E))$ is annihilated by some power of $t-\theta$. Hence in this case, $\dumot(E)$ is indeed an effective Anderson $t$-comotive.

\begin{rem}
We showed in \cite{am:aefam} for perfect fields $K$ that the conditions of being abelian and being coabelian are indeed equivalent. This result was extended in \cite{qg-am:pamfrfe} to more general base rings $R$ instead of the perfect field $K$.
\end{rem}

Conversely, assume we are given an effective Anderson $t$-motive $\mot$, and we encounter that it is also finitely generated free as $K\sp{\tau}$-module with some $K\sp{\tau}$-basis $(\bkappa_1,\ldots, \bkappa_d)$, and the $t$-action on $\mot$ can be described as 
\[  t\cdot \svect{\bkappa}{d} = D\cdot \svect{\bkappa}{d}.\]

Then this $t$-motive is isomorphic to the $t$-motive associated to the Anderson $t$-module $(E,\phi)$ with $E=\Ga^d$, and
\[  \phi_t = D \in \Mat_{d\times d}(K\sp{\tau})\cong \End_{\grp,\Fq}(\Ga^d). \]
Here the condition that $\mot(E)/\left(K[t]\cdot \tau(\mot(E))\right)$ is annihilated by some power of $t-\theta$ ensures that 
$\dphi_a-\ell(a)$ is nilpotent for all $a\in \Fq[t]$.

Similarly, but computationally a bit more technical, assume we are given an effective Anderson $t$-comotive $\dumot$ over a perfect field $K$, and we encounter that it is also finitely generated free as $K\sp{\sigma}$-module.
If we have computed a $K\sp{\sigma}$-basis and the $t$-action on that basis, it is straight forward to obtain (up to isomorphism) the Anderson $t$-module to which this $t$-comotive is associated to.

In order to do so, one just has to reverse the steps given at the beginning of Section \ref{sec:dual-t-motives}.

\section{Bases of abelian $t$-motives}\label{sec:basis}

In this section, we answer Question \ref{item:basis-of-t-motive} of the introduction.

Throughout this section, we fix a $t$-module $(E,\phi)$  over $K$ of dimension $d$, as well as a coordinate system $\bkappa$, i.e., an $\Fq$-linear isomorphism of group schemes $\bkappa:E\cong \Ga^d$ defined over $K$, and just write $\mot$ for $\mot(E)$.

Let $\bkappa_i:E\to \Ga$ be the $i$-th component of $\bkappa$, i.e., the composition of $\bkappa$ with the projection $\pr_i:\Ga^d\to \Ga$ to the $i$-th component of $\Ga^d$. 
Then the set $\{\bkappa_1,\ldots, \bkappa_d\}$ is a $K\{\tau\}$-basis of the $t$-motive of $\mot$.

Furthermore, with respect to this coordinate system, we can represent $\phi_t$ by a matrix $D=(D_{ij})\in \Mat_{d\times d}(K\{\tau\})$, and the $t$-action on $\mot$ is given in matrix notation by\footnote{See also \cite[Section 5]{am:aefam} for more details.}
\begin{equation} \label{eq:t-action-on-mot} 
t\cdot \svect{\bkappa}{d} = D\cdot \svect{\bkappa}{d}.
\end{equation}

So we are exactly in the situation of Section \ref{sec:finite-generation-and-relations} with $\fop=\tau$ and $\sop=t$, as well as $\gamma_\tau(x)=x^q$ for all $x\in K$, and $\gamma_t=\id_K$. We therefore view $\mot$ as the quotient $\cF/\Dgen{p_1,\ldots, p_d}$ with 
\[ \cD=K\sp{\tau,t},\quad \cF=\bigoplus_{i=1}^d \cD\kappa_i,\]
and for $i=1,\ldots, d$:
\[  p_i = t\kappa_i - \sum_{j=1}^d D_{ij}\kappa_j\in \cF. \]

From Theorem \ref{thm:criterion-for-finiteness} we obtain

\begin{thm}\label{thm:criterion-for-abelianess}
Let $K\sp{\tau,t}$ be equipped with the lexicographical order with $t\prec \tau$, and extend this monomial order to $\cF$ via a position-over-term order. Let $J$ be a Janet basis of $\Dgen{p_1,\ldots, p_d}$ obtained by applying Algorithm \ref{algo:janet-basis}.
Let the quantities $\sfn_i,\sfm_i$ ($i=1,\ldots,d$), $\Wgen$, $\bWgen$, $\Btop$ and $\Blow$ be given as in Definition \ref{def:quantities} and Theorem \ref{thm:relations-and-fop-action} for $\fop=\tau$ and $\sop=t$.

Then
\begin{enumerate}
\item \label{item:mot-abelian} $\mot$ is abelian $\Leftrightarrow$ $\forall i=1,\ldots, d$, we have $\sfn_i<\infty$.
\item \label{item:rational-t-motive}
The rational $t$-motive $K(t)\otimes_{K[t]} \mot$ is a finite dimensional $K(t)$-vector space of dimension
\[ \dim_{K(t)}(K(t)\otimes_{K[t]} \mot) = \sum_{i=1}^d \sfm_i. \]
\end{enumerate}
If $\mot$ is abelian, we further have:
\begin{enumerate} \setcounter{enumi}{2}
\item \label{item:rank} the rank of $\mot$ as a $K[t]$-module is given by
 \[ \rk_{K[t]}(\mot) = \sum_{i=1}^d \sfm_i.\footnote{In \cite[Definition 7.2]{am:naacigr}, we also defined the notion of a \emph{rank} for non-abelian $t$-motives, namely as the dimension of the rational $t$-motive $K(t)\otimes_{K[t]} \mot$. We nevertheless decided to add the condition ``abelian'' here, as the more general notion is not so common. } \]
\item \label{item:motive-relations} $\bWgen$ is a (finite) generating set for $\mot$ as $K[t]$-module, and their $K[t]$-relations are generated by the elements in $\Blow$ (written as $K[t]$-linear combinations of the elements in $\Wgen$).
\item \label{item:tau-action} The $\tau$-action is given as follows: For $w\in \Wgen$, i.e., $w=\tau^j\kappa_i$ for some $i=1,\ldots, d$ and $0\leq j< \sfn_i$, one has
\[   \tau(\tau^j\bkappa_i)= \partdef{ \tau^{j+1}\bkappa_i & \text{if } j<\sfn_i-1\\ \tau^{\sfn_i}\bkappa_i-\frac{1}{\lc(b_i)}\overline{b_i} & \text{if } j=\sfn_i-1\,  } \]
where $b_i\in \Btop$ is the element with leading monomial $\lm(b_i)=\tau^{\sfn_i}\kappa_i$, and $\overline{b_i}$ is obtained from it by replacing $\kappa_i$ with $\bkappa_i$.
\end{enumerate}
\end{thm}

\begin{proof}
Part \eqref{item:mot-abelian} follows directly from Theorem \ref{thm:criterion-for-finiteness}\eqref{item:M-finite} and \eqref{item:sop-rank} taking into account that $\gamma_t=\id_K$ is an isomorphism.

Parts \eqref{item:rational-t-motive} and \eqref{item:rank} are nothing else than Parts \eqref{item:rational-dimension} and \eqref{item:sop-rank} of Theorem \ref{thm:criterion-for-finiteness}.

Finally, the Parts \eqref{item:motive-relations} and \eqref{item:tau-action} are rephrasings of Theorem \ref{thm:relations-and-fop-action}.
\end{proof}

\begin{rem}
In \cite{am:aefam}, we also provided an algorithm to check whether a $t$-module is abelian which we shortly recall here. Assume that $K$ is perfect, and consider the entries of the matrix $t\cdot \one_d-D$ as elements of polynomials in $t$ over the skew Laurent series ring in $\sigma=\tau^{-1}$,
\[ K\sls{\sigma}=\biggl\{ \sum_{i=i_0}^\infty x_i \sigma^i \mid i_0\in \ZZ,x_i\in K \biggr\} \]
using the embedding
\[   K\sp{\tau} \to K\sls{\sigma}, \sum_{j=0}^n x_j \tau^j\mapsto \sum_{j=0}^n x_j \sigma^{-j}. \]
Diagonalize the matrix $t\cdot \one_d-D$ via row and column operations in such a way that the diagonal entries successively divide each other.
Then consider the last entry $\lambda_d$ of the resulting diagonal matrix $\diag(\lambda_1,\ldots, \lambda_d)$ which by assumption is divisible by all other $\lambda_i$ ($i=1,\ldots,d-1$). The $t$-module is abelian if and only if the Newton polygon of that polynomial $\lambda_d\in K\sls{\sigma}[t]$ has positive slopes only.

An interesting question is how that algorithm is related to the Janet algorithm.
The starting elements $p_1,\ldots, p_d$ for the Janet algorithm are the rows of the matrix $t\cdot \one_d-D$. As multiplication with $\tau$ is an isomorphism when considered in $K\sls{\sigma}[t]$, computing the normal forms of $\tau$-multiples of some $p_i$ look like computing the result of an elementary row operation. 
However, as other multiples of $p_i$ could occur in the reduction process, in general it is not.
Furthermore, there are no transformations equivalent to column operations in the Janet algorithm.

So at first sight, there are similarities, but also differences, and it would need a thorough investigation to clarify the situation.
\end{rem}

\section{Examples} \label{sec:examples}

We illustrate the algorithm and the results by the examples that were also considered in \cite[Section 9]{am:aefam}.

\begin{exmp}\label{exmp:drinfeld-module}
Let $\psi$ be a Drinfeld module over $K$ of rank $r$ and characteristic $\ell:\Fq[t]\to K$, i.e.,
\[ \psi_t=\theta+a_1\tau+\ldots + a_r\tau^r \]
with $\theta:=\ell(t)\in K$ and $a_1,\ldots, a_r\in K$, $a_r\ne 0$. We will see that in this case, everything becomes quite simple and fits perfectly to the well-known.

We have $d=1$, and
\[ \cF = \cD\kappa_1/\Dgen{p_1} \]
where $p_1=t\kappa_1 - \psi_t\cdot \kappa_1$.
As we only have one generator $p_1$, the Janet algorithm terminates directly producing $J=\{ (p_1,\mu_1) \}$ with $\mu_1=\Mon(\{\tau,t\})$.
The leading monomial of $p_1$ is $\tau^r\kappa_1$, hence
\[ \sfn_1 = \sfm_1 = r,\quad  \Wgen=\Wind = \{ \tau^j\kappa_1 \mid 0\leq j< r \}. \]
So, we have $r$ different $K[t]$-generators and no relations, i.e., $\{ \tau^j\bkappa_1 \mid 0\leq j< r \}$ is a $K[t]$-basis for the $t$-motive $\mot(E)$, and the $\tau$-action is given by
\[ \tau(\tau^j\bkappa_1) = \partdef{
\tau^{j+1}\bkappa_1, & j<r-1 \\
\pr\left(\tau^{r}\kappa_1-\frac{1}{a_r}p_1\right) & \\ 
\qquad = \frac{t-\theta}{a_r}\bkappa_1 - \frac{a_1}{a_r}\tau\bkappa_1-\ldots -\frac{a_{r-1}}{a_r}\tau^{r-1}\bkappa_1,\quad     & j=r-1.
} \]
\end{exmp}

\begin{exmp}\label{exmp:quasi-periodic-extension}
Quasi-periodic extensions of Drinfeld modules by $\Ga$:\\
Let $\psi$ be a Drinfeld module over $K$ of rank $r$ and characteristic $\ell:\Fq[t]\to K$, and let $\delta:\Fq[t]\to \tau K\{\tau\},a\mapsto \delta_a$ be a $\ell$-$\psi$-bi-derivation, i.e., an $\Fq$-linear map satisfying
\[ \delta_{ab}=\ell(a)\cdot \delta_b+\delta_a\cdot \psi_b \]
for all $a,b\in \Fq[t]$.
We consider the $t$-module $(\Ga^2,\phi)$ of dimension $2$ given by

\[ \phi_t = \begin{pmatrix} \psi_t & 0 \\ \delta_t & \theta \end{pmatrix},\]
where $\theta:=\ell(t)$.

Here, $d=2$, $p_1=t\kappa_1-\psi_t\kappa_1$ and $p_2=t\kappa_2-\delta_t\kappa_1-\theta\kappa_2$.
The leading monomials are $\lm(p_1)=-\lm(\psi_t)\kappa_1=\tau^r\kappa_1$ and $\lm(p_2)=-\lm(\delta_t)\kappa_1$.
So we would have to successively reduce one by the other until we get a reduced system.

Here, it is much more efficient to change the coordinate system by swapping $\kappa_1$ and $\kappa_2$ (cf.~Remark \ref{rem:changing-coordinate-system}). For an easier reading, we will not change the indices, but do the computation with $\kappa_2\succ \kappa_1$. Namely, in this ordering, we have still have
$\lm(p_1)=\tau^r\kappa_1$, but $\lm(p_2)=t\kappa_2$. So, if $\deg_{\tau}(\delta_t)<r$, we already have reduced polynomials, and if $\deg_{\tau}(\delta_t)\geq r$, we have to reduce $p_2$ by $p_1$ once to obtain reduced polynomials. In both cases the leading monomials don't change, and we directly obtain a Janet basis
\[ J= \{ (p_2,\{\tau,t\}), (p_1,\{\tau,t\}) \}, \]
-- depicted below --
from which we can read off that $\sfn_2=\infty$ (the $\sfn$ corresponding to $\kappa_2$), i.e., that this $t$-module is not abelian.

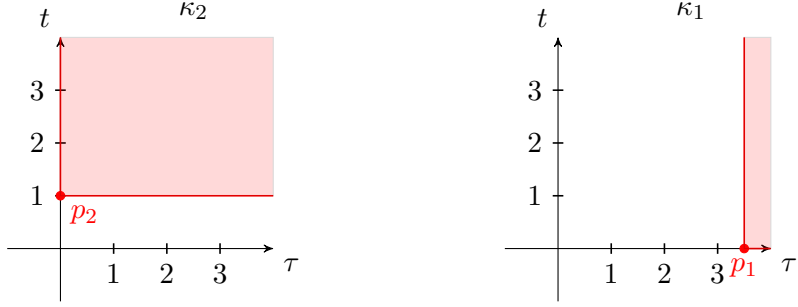
\twosheets{
	\grid{\tau}{t}{2}
	\point{p_2}{below right}{0}{1}
    \cone{0}{1} 
}{
	\grid{\tau}{t}{1}
	\point{p_1}{below}{3.5}{0}
    \cone{3.5}{0}
}

Although, the $t$-module is not abelian, the Janet basis gives us a description of the $t$-motive by generators and relations. Namely from Proposition \ref{prop:properties-of-quantities}, we know that $\mot$ is generated by 
\[ \bWgen = \{ \tau^j\kappa_2 \mid j\geq 0\}\cup \{ \tau^j\kappa_1 \mid 0\leq j<r=\sfn_1 \}. \]
The $K[t]$-relations are not only generated by $p_2$, but here also by all $\tau$-shifts of $p_2$. Hence,
\[ (t-\theta)\bkappa_2 = \delta_t\bkappa_1, \]
and for all $j\geq 1$,
\[ (t-\theta^{q^j})\tau^j\bkappa_2 =\tau^j\delta_t\bkappa_1=\NF(\tau^j\delta_t\bkappa_1, \{ p_1\})\in \Ktgen{\bkappa_1,\ldots, \tau^{r-1}\bkappa_1}. \]
This verifies Theorem \ref{thm:criterion-for-abelianess}\eqref{item:rational-t-motive}, that $K(t)\otimes_{K[t]} \mot$ is finitely generated of dimension $r$.
\end{exmp}

\begin{exmp}\label{exmp:carlitz-tensor-power}

The second Carlitz tensor power $C^{\otimes 2}$ is given by 
\[ D=\begin{pmatrix} \theta & 1 \\ \tau & \theta \end{pmatrix}. \]
So we have $d=2$, and
\[ p_1 = t\kappa_1 - \theta\kappa_1 -\kappa_2,\quad p_2=-\tau\kappa_1+t\kappa_2-\theta\kappa_2. \]
\begin{itemize}
\item First loop in Janet algorithm: 
\begin{itemize}
\item $G_1=\{p_1,p_2\}$ already reduced. 
\item Janet decomposition $J_1= \{ (p_2,\{\tau,t\}), (p_1,\{t\}) \}$,
\item $P_1= \{ \NF(\tau p_1, J_1)\} = \{ -\tau\kappa_2+ (t-\theta)(t-\theta^q)\kappa_2 \}$,
\end{itemize}
\item Second loop in Janet algorithm:
\begin{itemize}
\item $G_2=\{p_1,p_2,p_3\}$ with $p_3=-\tau\kappa_2+ (t-\theta)(t-\theta^q)\kappa_2$,
\item $J_2= \{ (p_2,\{\tau,t\}), (p_1,\{t\}), (p_3,\{\tau,t\})\}$,
\item $P_2= \{ \NF(\tau p_1, J_2)\}\setminus \{0\} = \emptyset$.
\end{itemize}
\end{itemize}
Hence, $J= \{ (p_2,\{\tau,t\}), (p_1,\{t\}), (p_3,\{\tau,t\})\}$ is a Janet basis.

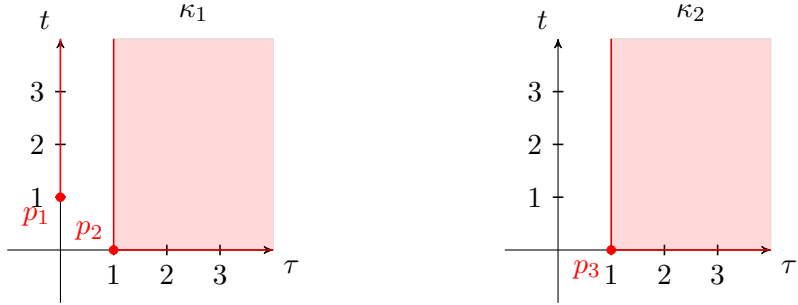
\twosheets{
	\grid{\tau}{t}{1}
	\point{p_2}{above left}{1}{0}
	\point{p_1}{below left}{0}{1}
    \linet{0}{1} \cone{1}{0}
}{
	\grid{\tau}{t}{2}
	\point{p_3}{below left}{1}{0}
    \cone{1}{0}
}
 
We obtain $\sfn_1=\sfn_2=1$, and $\sfm_1=0$, $\sfm_2=1$. Hence, $\mot=\mot(C^{\otimes 2})$ is abelian of rank $\sfm_1+\sfm_2=1$.
\begin{itemize}
\item $B:=\{ \bkappa_1, \bkappa_2\}$ generates $\mot$ as $K[t]$-module, 
\item the relation is given by $0=p_1=(t-\theta)\bkappa_1-\bkappa_2$, i.e., $\bkappa_2=(t-\theta)\bkappa_1$, and hence $\{\bkappa_1\}$ is a basis.
\item the $\tau$-action is given by
\[ \tau(\bkappa_1) = \tau\bkappa_1+\overline{p_2}=(t-\theta)\bkappa_2=(t-\theta)^2\bkappa_1.\]
\end{itemize}

\medskip

A similar calculation works for the $d$-th Carlitz tensor power $C^{\otimes d}$ which is given by
\[ D =\begin{pmatrix} 
\theta & 1 & 0 & \cdots & 0\\ 
0 & \ddots & \ddots & \ddots & \vdots\\
\vdots & \ddots &  \ddots & \ddots & 0 \\
0 & & \ddots & \ddots & 1 \\
\tau & 0 & \cdots & 0 & \theta \end{pmatrix} \in \Mat_{d\times d}(K\{\tau\})
.\]
Of course, the calculation reveals the well-known result that its $t$-motive has $\kappa_1$ as its $K[t]$-basis, and $\tau(\kappa_1)=(t-\theta)^d\kappa_1$.
We leave the detailed calculation to the reader.
\end{exmp}

\begin{exmp}\label{exmp:special-example}
We consider Example 6.3 from \cite{am:aefam}, i.e., the $t$-module $(E,\phi)$ over the rational function field $K=\Fq(\theta)$ with
\[ D= \begin{pmatrix} \theta+\tau^2 & \tau^3 \\ 1+\tau & \theta+\tau^2 \end{pmatrix}
 \]
with respect to a coordinate system $\{\bkappa_1,\bkappa_2\}$.
In this case,
\[  p_1=t\kappa_1-\tau^2\kappa_1-\theta\kappa_1 - \tau^3\kappa_2\quad \text{and}\quad p_2=t\kappa_2-\tau\kappa_1-\kappa_1-\tau^2\kappa_2-\theta\kappa_2. \]

As in the previous examples, the calculations become easier when we change the order of the basis to $\kappa_1\prec \kappa_2$.\footnote{We leave it to the reader to run the algorithm with the order $\kappa_1\succ \kappa_2$.}

\begin{itemize}
\item First loop in Janet algorithm: 
\begin{itemize}
\item $\lm(p_2)=\tau^2\kappa_2$ divides $\lm(p_1)=\tau^3\kappa_2$, so we have to reduce $p_1$, and obtain
\[  p_1'=p_1-\tau p_2= -\tau t\kappa_2+\theta^q\tau\kappa_2 + \tau \kappa_1+t\kappa_1-\theta \kappa_1, \]
and $G_1=\{p_2,p_1'\}$.
\item Janet decomposition $J_1= \{ (p_2,\{\tau,t\}), (p'_1,\{t\}) \}$, 
\item $p_3:=\NF(\tau p'_1, J_1) 
=-(t-\theta)(t-\theta^{q^2})\kappa_2+ (\tau^2+2\tau t-\theta^q\tau -\theta^{q^2}\tau +t-\theta^{q^2})\kappa_1$, hence $P_1=\{p_3\}$.
\end{itemize}
\item Second loop in Janet algorithm: 
\begin{itemize}
\item $G_2=\{p_2,p'_1,p_3\}$ is already auto-reduced, and
\item $ J_2= \{ (p_2,\{\tau,t\}), (p'_1,\{t\}), (p_3,\{t\})\}$,
\item $\NF(\tau p'_1, J_2)=0$, and
\[ p_4:=\NF(\tau p_3, J_2)= \left(\tau^3+(2t-\theta^{q^2}-\theta^{q^3})\tau^2 -(t-\theta)(t-\theta^{q^3})\right)\kappa_1, \]
hence $P_2=\{p_4\}$.
\end{itemize}
\item Third loop in Janet algorithm: 
\begin{itemize}
\item $G_3=\{p_2,p'_1,p_3,p_4\}$ is already auto-reduced, and
\item $ J_3= \{ (p_2,\{\tau,t\}), (p'_1,\{t\}), (p_3,\{t\}), (p_4,\{\tau,t\}) \}$,
\item $\NF(\tau p'_1, J_3)=0$, and $\NF(\tau p_3, J_3)=0$, hence $P_3=\emptyset$.
\end{itemize}
\end{itemize}
Therefore, $J= \{ (p_2,\{\tau,t\}), (p'_1,\{t\}), (p_3,\{t\}), (p_4,\{\tau,t\})\}$ is a Janet basis.

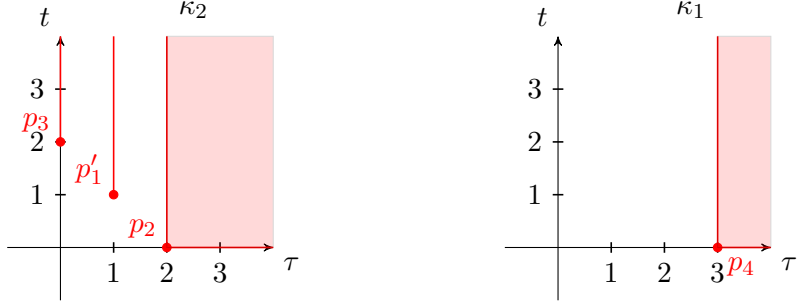
\twosheets{
	\grid{\tau}{t}{2}
	\point{p_2}{above left}{2}{0}
	\point{p'_1}{above left}{1}{1}
	\point{p_3}{above left}{0}{2}
    \cone{2}{0}
    \linet{1}{1} 
    \linet{0}{2} 
}{
	\grid{\tau}{t}{1}
	\point{p_4}{below right}{3}{0}
    \cone{3}{0}
}

The result shows that the $t$-motive is abelian of rank $3$, and
the set 
$\{ \tau \bkappa_2, \bkappa_2,\tau^2\bkappa_1,\tau \bkappa_1,\bkappa_1 \}$
is generating the $t$-motive. The relations are provided by $p'_1$ and $p_3$, i.e.,
\begin{align*}
\tau \bkappa_1 &= (t-\theta^q)\tau\bkappa_2 - (t-\theta) \bkappa_1,\quad \text{and} \\
\tau^2\bkappa_1 &= (t-\theta)(t-\theta^{q^2})\bkappa_2 - (2t-\theta^q -\theta^{q^2})\tau\bkappa_1-(t-\theta^{q^2})\bkappa_1 \\
&= (t-\theta)(t-\theta^{q^2})\bkappa_2 - (2t-\theta^q -\theta^{q^2})\left(  
(t-\theta^q)\tau\bkappa_2 - (t-\theta) \bkappa_1\right)\\
&\qquad -(t-\theta^{q^2})\bkappa_1 \\
&= (t-\theta)(t-\theta^{q^2})\bkappa_2 - (2t-\theta^q -\theta^{q^2})(t-\theta^q)\tau\bkappa_2 \\
&\qquad +\left( (t-\theta)(2t-\theta^q -\theta^{q^2}) - (t-\theta^{q^2})\right) \bkappa_1
\end{align*}
So, $\be_1=\tau\bkappa_2$, $\be_2=\bkappa_2$, $\be_3=\bkappa_1$ is a $K[t]$-basis. We compute the $\tau$-action using Theorem \ref{thm:criterion-for-abelianess}\eqref{item:tau-action} and eliminating $\tau \bkappa_1$ and $\tau^2\bkappa_1$ afterwards, and obtain
\begin{align*}
\tau(\be_1)&= \tau^2\bkappa_2 = -(t-\theta^q)\be_1+(t-\theta)\be_2+(t-\theta-1)\be_3, \\
\tau(\be_2)&= \be_1 \\
\tau(\be_3)&= (t-\theta^q)\be_1-(t-\theta)\be_3.
\end{align*}

\end{exmp}

\section{Computing the Anderson $t$-module to a $t$-motive}\label{sec:reverse-direction}

In Section \ref{sec:basis}, we started with an Anderson $t$-module from which we readily had the $t$-motive $\mot$ as a $K\sp{\tau}$-module with $t$-action, and computed a presentation of $\mot$ as a $K[t]$-module with a $\tau$-action.

In this section, we go the reverse direction, and answer Question \ref{item:t-module-to-t-motive} of the introduction.

\medskip

Let $(\mot,\tau_\mot)$ be an effective Anderson $t$-motive of rank $r$ as defined in Section \ref{sec:t-modules}. Let $\{\be_1,\ldots, \be_r\}$ be a $K[t]$-basis of $\mot$, and (written in matrix form)
\[ \tau_\mot \svect{\be}{r} = \Theta \svect{\be}{r} \]
for some $\Theta\in \Mat_{r\times r}(K[t])$.

So we are exactly in the situation of Section \ref{sec:finite-generation-and-relations}, this time, however, with $\fop=t$, and $\sop=\tau$, as well as $\be_1,\ldots, \be_r$ in place of $\bkappa_1,\ldots,\bkappa_d$.

We therefore view $\mot$ as the quotient $\cF/\Dgen{p_1,\ldots, p_r}$ with 
\[ \cD=K\sp{t,\tau},\quad \cF=\bigoplus_{i=1}^r \cD e_i,\]
and for $i=1,\ldots, r$:
\[  p_i = \tau e_i - \sum_{j=1}^r \Theta_{ij}e_j\in \cF. \]

\begin{thm}\label{thm:criterion-for-t-module-to-motive}
Let $K\sp{t,\tau}$ be equipped with the lexicographical order with $\tau\prec t$, and extend this monomial order to $\cF$ via a position-over-term order. Let $J$ be a Janet basis of $\Dgen{p_1,\ldots, p_r}$ obtained by applying Algorithm \ref{algo:janet-basis}.
Let the quantities $\sfn_i,\sfm_i$ ($i=1,\ldots,r$), $\Wgen$, $\bWgen$, $\Btop$ and $\Blow$ be given as in Definition \ref{def:quantities} and Theorem \ref{thm:relations-and-fop-action} for $\fop=t$ and $\sop=\tau$, as well as $\be_1,\ldots, \be_r$ in place of $\bkappa_1,\ldots,\bkappa_d$.

Then
\begin{enumerate}
\item \label{item:mot-tau-finite} $\mot$ is finitely generated as $K\sp{\tau}$-module $\Leftrightarrow$ $\forall i=1,\ldots, r$, we have $\sfn_i<\infty$.
\end{enumerate}
If $\mot$ is finitely generated as $K\sp{\tau}$-module, we further have
\begin{enumerate} \setcounter{enumi}{1}
\item \label{item:mot-tau-relations} $\bWgen$ is a (finite) generating set for $\mot$ as $K\sp{\tau}$-module, and the $K\sp{\tau}$-relations are generated by the elements in $\Blow$ (written as $K\sp{\tau}$-linear combinations of the elements in $\Wgen$).
\item \label{item:t-action} The $t$-action is given as follows: For $w\in \Wgen$, i.e., $w=t^j e_i$ for some $i=1,\ldots, r$ and $0\leq j< \sfn_i$, one has
\[   t(t^j\be_i)= \partdef{ t^{j+1}\be_i & \text{if } j<\sfn_i-1\\ t^{\sfn_i}\be_i-\frac{1}{\lc(b_i)}\overline{b_i} & \text{if } j=\sfn_i-1\,  } \]
where $b_i\in \Btop$ is the element with leading monomial $\lm(b_i)=t^{\sfn_i}e_i$, and $\overline{b_i}$ is obtained from it by replacing $e_i$ with $\be_i$.
\item \label{item:K-perfect-mot-free} If $K$ is perfect, then $\mot$ is free as $K\sp{\tau}$-module.
\end{enumerate}
If $\mot$ is finitely generated free as $K\sp{\tau}$-module, then
\begin{enumerate} \setcounter{enumi}{4}
\item \label{item:mot-of-E} $\mot$ is the $t$-motive associated to some (abelian) Anderson $t$-module $(E,\phi)$.
\item \label{item:dimension-of-E} The dimension of $E$ is given as
 \[ \dim(E) = \sum_{i=1}^r \sfm_i. \]
\end{enumerate}
\end{thm}

\begin{proof}
Parts \eqref{item:mot-tau-finite}--\eqref{item:K-perfect-mot-free} directly follow from Theorem \ref{thm:criterion-for-finiteness} and Theorem \ref{thm:relations-and-fop-action}.
That \eqref{item:mot-of-E} holds was explained in Section \ref{subsec:relation-between-Andersons-objects}, and \eqref{item:dimension-of-E} is  the rank formula in Theorem \ref{thm:criterion-for-finiteness}\eqref{item:sop-rank} taking into account that $\dim(E)=\rk_{K\sp{\tau}}(\mot)$.
\end{proof}

As explained in Section \ref{subsec:relation-between-Andersons-objects}, we can easily describe the Anderson $t$-module if we have found a $K\sp{\tau}$-basis of $\mot$ and the $t$-action on it.
Further, in Remark \ref{rem:finding-a-basis}, we have explained that obtaining a basis from a finite generating set is just the non-commutative version of the algorithm for the elementary divisor theorem.

\begin{rem}
In \cite{uh:iaamm}, a more general definition of an Anderson $t$-module $E$ is used, namely that $E$ does not have to be isomorphic to $\Ga^d$ itself, but have to become isomorphic after base change to some finite extension $L\supset K$.

If one uses that definition, then Theorem \ref{thm:criterion-for-t-module-to-motive}\eqref{item:mot-of-E} and \eqref{item:dimension-of-E} hold without the freeness assumption. Namely, Part \eqref{item:K-perfect-mot-free} ensures that $\mot$ becomes free over the perfection of $K$, and since $\mot$ is finitely generated, this already holds over some finite inseparable extension $L$ of $K$.
So from the data, we obtain some $E$ that becomes isomorphic to $\Ga^d$ over $L$.
\end{rem}

\section{Anderson comotives}\label{sec:dual-t-motives}

In this section, we answer Questions \eqref{item:basis-of-dual-t-motive} and \eqref{item:t-module-to-dual-t-motive} regarding the $t$-comotives.

The story runs parallel to the one for $t$-motives with small technicalities due to turning the right-$K\sp{\tau}$-action on $\dumot$ to a left-$K\sp{\sigma}$-action.

\subsection{Answer to Question \ref{item:basis-of-dual-t-motive}}

Throughout this section, $K$ is perfect, and we fix an Anderson $t$-module $(E,\phi)$  over $K$ of dimension $d$, as well as a coordinate system $\bkappa$, i.e., an $\Fq$-linear isomorphism of group schemes $\bkappa:E\cong \Ga^d$ defined over $K$, and just write $\dumot$ for $\dumot(E)$.

Let $\bdk_j:\Ga\to E$ be the composition of the injection $\inj_j:\Ga\to \Ga^d$ into the $j$-th component with $\bkappa^{-1}$ for all $j=1,\ldots, d$.
Then the tuple $(\bdk_1,\ldots, \bdk_d)$ is a $K\{\sigma\}$-basis of the $t$-comotive of $\dumot$.

Furthermore, with respect to this coordinate system, we can represent $\phi_t$ by a matrix $D=(D_{ij})\in \Mat_{d\times d}(K\{\tau\})$, and the $t$-action on $\dumot$ is described by 
\begin{equation}\label{eq:t-action-on-dumot} 
 t\cdot \zvect{\bdk}{d} = \zvect{\bdk}{d}\cdot D
\end{equation}
in terms of the right-$K\sp{\tau}$-action, and hence in terms of the left-$K\sp{\sigma}$-action as
\begin{equation}\label{eq:t-action-on-dumot-2} 
 t\cdot \svect{\bdk}{d} = D^* \svect{\bdk}{d}.
\end{equation}
Here $D^*$ is obtained from $D$ by transposing and by transforming the entries $f=\sum_{i} a_i\tau^i\in K\sp{\tau}$ to 
$f^*= \sum_{i} \sigma^i a_i=\sum_{i} (a_i)^{1/q^i}\sigma^i\in K\sp{\sigma}$.

Again, we are exactly in the situation of Section \ref{sec:finite-generation-and-relations}, this time with $\fop=\sigma$ and $\sop=t$, as well as $\gamma_\sigma(x)=x^{1/q}$ for all $x\in K$, and $\gamma_t=\id_K$. We therefore view $\dumot$ as the quotient $\cF/\Dgen{p_1,\ldots, p_d}$ with 
\[ \cD=K\sp{\sigma,t},\quad \cF=\bigoplus_{i=1}^d \cD\dk_i,\]
and for $i=1,\ldots, d$:
\[  p_i = t\dk_i - \sum_{j=1}^d D^*_{ij}\dk_j\in \cF. \]

\begin{thm}\label{thm:criterion-for-t-finiteness}
Let $K\sp{\sigma,t}$ be equipped with the lexicographical order with $t\prec \sigma$, and extend this monomial order to $\cF$ via a position-over-term order. Let $J$ be a Janet basis of $\Dgen{p_1,\ldots, p_d}$ obtained by applying Algorithm \ref{algo:janet-basis}.
Let the quantities $\sfn_i,\sfm_i$ ($i=1,\ldots,d$), $\Wgen$, $\bWgen$, $\Btop$ and $\Blow$ be given as in Definition \ref{def:quantities} and Theorem \ref{thm:relations-and-fop-action} for $\fop=\sigma$, $\sop=t$, and $\dk_1,\ldots,\dk_d$ in place of $\kappa_1,\ldots,\kappa_d$.

Then
\begin{enumerate}
\item \label{item:dumot-t-finite} $\dumot$ is coabelian $\Leftrightarrow$ $\forall i=1,\ldots, d$, we have $\sfn_i<\infty$.
\item \label{item:rational-dual-t-motive}
The rational $t$-comotive $K(t)\otimes_{K[t]} \dumot$ is a finite dimensional $K(t)$-vector space of dimension
\[ \dim_{K(t)}(K(t)\otimes_{K[t]} \dumot) = \sum_{i=1}^d \sfm_i. \]
\end{enumerate}
If $\dumot$ is coabelian, then
\begin{enumerate} \setcounter{enumi}{2}
\item \label{item:rank-dumot} its rank as a $K[t]$-module is given by
 \[ \rk_{K[t]}(\dumot) = \sum_{i=1}^d \sfm_i.\footnote{In \cite[Definition 7.2]{am:naacigr}, we also defined the notion of a \emph{rank} for non-coabelian $t$-comotives, namely as the dimension of the rational $t$-comotive $K(t)\otimes_{K[t]} \dumot$. We nevertheless decided to add the condition ``coabelian'' here, as the more general notion is not so common. } \]
\item \label{item:dual-motive-relations} $\bWgen$ is a (finite) generating set for $\dumot$ as $K[t]$-module, and the $K[t]$-relations are generated by the elements in $\Blow$ (written as $K[t]$-linear combinations of the elements in $\Wgen$).
\item \label{item:sigma-action} The $\sigma$-action is given as follows: For $w\in \Wgen$, i.e., $w=\sigma^j\dk_i$ for some $i=1,\ldots, d$ and $0\leq j< \sfn_i$, one has
\[   \sigma(\sigma^j\bdk_i)= \partdef{ \sigma^{j+1}\bdk_i & \text{if } j<\sfn_i-1\\ \sigma^{\sfn_i}\bdk_i-\frac{1}{\lc(b_i)}\overline{b_i} & \text{if } j=\sfn_i-1\,  } \]
where $b_i\in \Btop$ is the element with leading monomial $\lm(b_i)=\sigma^{\sfn_i}\dk_i$, and $\overline{b_i}$ is obtained from it by replacing $\dk_i$ with $\bdk_i$.
\end{enumerate}
\end{thm}

\begin{proof}
The proof is the same as for Theorem \eqref{thm:criterion-for-abelianess}.
\end{proof}

\begin{exmp}\label{exmp:special-example-dual-t-motive}
We compute a $K[t]$-basis of the $t$-comotive associated to the Anderson $t$-module $E$ of Example \ref{exmp:special-example}, i.e., where
\[ D= \begin{pmatrix} \theta+\tau^2 & \tau^3 \\ 1+\tau & \theta+\tau^2 \end{pmatrix}
 \]
with respect to a coordinate system $\kappa$.
So with respect to the corresponding basis  $\{\dk_1,\dk_2\}$, the $t$-action is
\[ t \cdot \begin{pmatrix}\dk_1 \\ \dk_2\end{pmatrix} = D^*\cdot \begin{pmatrix}\dk_1 \\ \dk_2\end{pmatrix} = \begin{pmatrix} \theta+\sigma^2 & 1+\sigma \\ \sigma^3 & \theta+\sigma^2 \end{pmatrix} \begin{pmatrix}\dk_1 \\ \dk_2\end{pmatrix}. \]
Therefore,
\[  p_1=t\dk_1-\sigma^2\dk_1-\theta\dk_1 - \sigma\dk_2-\dk_2 \quad \text{and}\quad p_2=t\dk_2-\sigma^3\dk_1-\sigma^2\dk_2-\theta\dk_2. \]

Comparing this with the situation in Example \ref{exmp:special-example}, we see that the dictionary
\[  \kappa_1\leftrightarrow \dk_2,\quad \kappa_2\leftrightarrow \dk_1,\quad \tau\leftrightarrow \sigma \]
maps the relation $p_1$ here to the relation $p_2$ there, and the relation $p_2$ here to $p_1$ there.

Hence, by choosing the order $\dk_1\succ \dk_2$, we obtain the very same computations in the Janet algorithm, only with $q$-power roots of $\theta$ instead of $q$-powers of $\theta$. So we end up with the  Janet basis
\[ J= \{ (p_1,\{\sigma,t\}), (p'_2,\{t\}), (p_3,\{t\}), (p_4,\{\sigma,t\})\},\]
where
\begin{align*}
 p_1 &= -\sigma^2\dk_1+(t-\theta)\dk_1-(1+\sigma)\dk_2,\\
 p'_2&= -\sigma t\dk_1+\theta^{1/q}\sigma \dk_1 + \sigma \dk_2+ (t-\theta)\dk_2,\\
 p_3 &= -(t-\theta)(t-\theta^{1/q^2})\dk_1+\left(\sigma^2+2\sigma t -\theta^{1/q}\sigma-\theta^{1/q^2}\sigma + t- \theta^{1/q^2}\right)\dk_2,\\
 p_4 &= \left(\sigma^3+(2t-\theta^{1/q^2}-\theta^{1/q^3})\sigma^2-(t-\theta)(t-\theta^{1/q^3})\right) \dk_2.
\end{align*}

So $\dumot$ is coabelian of rank $3$, and the set 
$\{ \sigma\dk_1, \dk_1, \sigma^2\dk_2,\sigma\dk_2,\dk_2\}$ is generating the $t$-comotive as $K[t]$-module. The relations are provided by $p'_2$ and $p_3$, showing that $\ce_1= \sigma\dk_1,\ce_2=\dk_1,\ce_3=\dk_2$ is a $K[t]$-basis.
The $\sigma$-action is given by 
\begin{align*}
\sigma(\ce_1) &= \sigma^2 \dk_1 = -(t-\theta^{1/q})\ce_1 + (t-\theta)\ce_2 + (t-\theta-1)\ce_3,\\
\sigma(\ce_2) &= \ce_1,\\
\sigma(\ce_3) &= (t-\theta^{1/q})\ce_1 - (t-\theta)\ce_3.
\end{align*}
\end{exmp}

\subsection{Answer to Question \ref{item:t-module-to-dual-t-motive}}

Let $(\dumot,\sigma_\dumot)$ be an effective Anderson $t$-comotive of rank $r$ as defined in Section \ref{sec:t-modules}. Let $\{\be_1,\ldots, \be_r\}$ be a $K[t]$-basis of $\dumot$, and (written in matrix form)
\[ \sigma_\dumot \svect{\be}{r} = \Theta \svect{\be}{r} \]
for some $\Theta\in \Mat_{r\times r}(K[t])$.

So we are exactly in the situation of Section \ref{sec:finite-generation-and-relations}, this time, however, with $\fop=t$, and $\sop=\sigma$, as well as $\be_1,\ldots, \be_r$ in place of $\bkappa_1,\ldots,\bkappa_d$.

We therefore view $\dumot$ as the quotient $\cF/\Dgen{p_1,\ldots, p_r}$ with 
\[ \cD=K\sp{t,\sigma},\quad \cF=\bigoplus_{i=1}^r \cD e_i,\]
and for $i=1,\ldots, r$:
\[  p_i = \sigma e_i - \sum_{j=1}^r \Theta_{ij}e_j\in \cF. \]

\begin{thm}\label{thm:criterion-for-t-module-to-dual-motive}
Let $K\sp{t,\sigma}$ be equipped with the lexicographical order with $\sigma\prec t$, and extend this monomial order to $\cF$ via a position-over-term order. Let $J$ be a Janet basis of $\Dgen{p_1,\ldots, p_r}$ obtained by applying Algorithm \ref{algo:janet-basis}.
Let the quantities $\sfn_i,\sfm_i$ ($i=1,\ldots,r$), $\Wgen$, $\bWgen$, $\Btop$ and $\Blow$ be given as in Definition \ref{def:quantities} and Theorem \ref{thm:relations-and-fop-action} for $\fop=t$ and $\sop=\sigma$, as well as $\be_1,\ldots, \be_r$ in place of $\bkappa_1,\ldots,\bkappa_d$.

Then
\begin{enumerate}
\item \label{item:dumot-sigma-finite} $\dumot$ is finitely generated as $K\sp{\sigma}$-module $\Leftrightarrow$ $\forall i=1,\ldots, r$, we have $\sfn_i<\infty$.
\end{enumerate}
If $\dumot$ is finitely generated as $K\sp{\sigma}$-module, we further have
\begin{enumerate} \setcounter{enumi}{1}
\item \label{item:dumot-sigma-relations} $\bWgen$ is a (finite) generating set for $\dumot$ as $K\sp{\sigma}$-module, and the $K\sp{\sigma}$-relations are generated by the elements in $\Blow$ (written as $K\sp{\sigma}$-linear combinations of the elements in $\Wgen$).
\item \label{item:t-action-on-dumot} The $t$-action is given as follows: For $w\in \Wgen$, i.e., $w=t^j e_i$ for some $i=1,\ldots, r$ and $0\leq j< \sfn_i$, one has
\[   t(t^j\be_i)= \partdef{ t^{j+1}\be_i & \text{if } j<\sfn_i-1\\ t^{\sfn_i}\be_i-\frac{1}{\lc(b_i)}\overline{b_i} & \text{if } j=\sfn_i-1\,  } \]
where $b_i\in \Btop$ is the element with leading monomial $\lm(b_i)=t^{\sfn_i}e_i$, and $\overline{b_i}$ is obtained from it by replacing $e_i$ with $\be_i$.
\item \label{item:dumot-free} $\dumot$ is free as $K\sp{\sigma}$-module.
\item \label{item:dumot-of-E} $\dumot$ is the $t$-comotive associated to some (coabelian) Anderson $t$-module $(E,\phi)$.
\item \label{item:dimension-of-E-for-dumot} The dimension of $E$ is given as
 \[ \dim(E) = \sum_{i=1}^r \sfm_i. \]
\end{enumerate}
\end{thm}

\begin{proof}
The proof is the same as the one for Theorem \ref{thm:criterion-for-t-module-to-motive}. The only difference is that $\gamma_\sigma$ is an isomorphism from the beginning, so Proposition \ref{prop:observations-on-M} already ensures that $\dumot$ is free as $K\sp{\sigma}$-module if it is finitely generated.
\end{proof}

As mentioned in Section \ref{subsec:relation-between-Andersons-objects}, we can easily describe the Anderson $t$-module if we have found a $K\sp{\sigma}$-basis of $\dumot$ and the $t$-action on it. This is done by reversing the steps at the beginning of this section where the $t$-comotive associated to the Anderson $t$-module was computed.

Further, in Remark \ref{rem:finding-a-basis}, we have explained that obtaining a basis from a finite generating set is just the non-commutative version of the algorithm for the elementary divisor theorem.

\subsection*{Acknowledgements}
The author would like to thank the referee for his suggestions for improving the presentation of this article.

\subsection* {Funding}
The author has received support from the SFB-TRR 195 ``Symbolic Tools in Mathematics and their Application'' of the German Research Foundation (DFG).


\normalsize

\def\cprime{$'$}

\end{document}